\documentclass[review]{elsarticle}

\usepackage{amssymb}
\setcitestyle{numbers,super}
\usepackage{amsmath}
\usepackage{amssymb}
\usepackage{amsfonts}
\usepackage{algorithmic}
\usepackage{graphicx}
\usepackage{textcomp}
\usepackage{xcolor}
\usepackage{graphicx}
\usepackage{subcaption}
\usepackage{mwe}
\usepackage{color, colortbl}
\usepackage{url}

\begin{document}
\begin{frontmatter}

\title{Risk-Constrained Community Battery Utilisation Optimisation for Electric Vehicle Charging with Photovoltaic Resources}


\author[inst1]{Khalil Gholami}
\ead{k.gholami@deakin.edu.au}
\affiliation[inst1]{organization={School of Eng, Deakin University},
            addressline={75 Pigdons Rd}, 
            city={Geelong},
            postcode={3216}, 
            state={Victoria},
            country={Australia}}

\author[inst2]{Asef Nazari}
\ead{asef.nazari@deakin.edu.au}
\affiliation[inst2]{organization={School of IT, Deakin University},
            addressline={75 Pigdons Rd}, 
            city={Geelong},
            postcode={3216}, 
            state={Victoria},
            country={Australia}}
\author[inst2]{Dhananjay Thiruvady}
\ead{dhananjay.thiruvady@deakin.edu.au}
\author[inst2]{Valeh Moghaddam}
\ead{valeh.moghaddam@deakin.edu.au}
\author[inst2]{Sutharshan Rajasegarar}
\ead{sutharshan.rajasegarar@deakin.edu.au}
\author[inst2]{Wei-Yu Chiu}
\ead{weiyu.chiu@deakin.edu.au}

\begin{abstract}
High penetration of renewable energy generation in the electricity grid presents power system operators with challenges including voltage instability mainly due to fluctuating power generation. To cope with intermittent renewable generation, community batteries introduce an elegant solution for storing excess generation of renewable resources and reverting to the grid in peak demand periods. The question of the right battery size coupled with the right investment Furthermore, the growth in adapting electrical vehicles (EVs) imposes additional demand-related challenges on the power system compared to traditional industrial and household demand. This paper introduces long-term planning for community batteries to capture the surplus generation of rooftop PV resources for a given area and redirect these resources to charge EVs, without direct injection to the upstream grid. For long-term investment planning on batteries, we consider 15 years' worth of historical data associated with solar irradiance, temperature, EV demands, and household demands. A novel stochastic mathematical model is proposed for decision-making on battery specifications (the type and capacity of battery per year) based on the four standard battery types provided by the Commonwealth Scientific and Industrial Research Organisation (CSIRO) in Australia. Uncertainties related to the EVs and RESs are captured by a non-parametric robust technique, named information gap decision theory, from optimistic and pessimistic perspectives. The investment decision-making part is formulated as mixed-integer linear programming taking advantage of the powerful commercial solver -- GUROBI -- which leads to finding feasible global solutions with low computational burden. The outcomes of this investigation not only detect optimal battery installation strategies to improve the stability profile of the grid by capturing the excess generation of PV resources but also facilitate EV integration in the community toward reaching net-zero emissions targets.
\end{abstract}


\end{frontmatter}

\bibliographystyle{elsarticle-num}

\section{Introduction}
Australia is recognised globally as a leader in solar energy adoption, securing the highest per-capita installation of solar panels \cite{riley2023connected}. Nearly one-third of Australian homes have solar panels to harness sun radiation to generate clean, renewable energy. As the technology continues to be cost-effective, solar panels are expected to become an even more prominent feature on Australian rooftops photovoltaic (PV) in the residential districts, further contributing to the country's clean energy goals \cite{li2020review}. This trend has grown so appealing in recent years that the use of PV resources has increased significantly in existing residential areas to decrease their electricity bills as well as to achieve a cleaner environment \cite{malekpour2016dynamic,murray2021voltage}. 

While renewable energy resources like solar, wind, and geothermal positively contribute to the sustainable future of every country, their large-scale integration into the power grid has many side effects \cite{gholami2022fuzzy,moradi2015optimal,dong2020distorted}. The variability and intermittency in renewable generation, as a function of weather conditions, make the forecasting of the generation hard. The lack of perfect control of the renewably generated energy causes imbalances in the electricity transmission and distribution systems affecting the frequency and voltage stability \cite{rubanenko2020analysis}. The stability of voltage and frequency is a crucial aspect of any power system as it impacts all other important characteristics of a power system including the reliability, efficiency, and safety of delivering electricity to customers. 

Electric vehicles (EVs) are already widely adopted, revolutionising the transportation industry by drastically lowering greenhouse gas emissions and reliance on fossil fuels \cite{ali2021maximizing,jordehi2020energy}. However, the increasing adoption of EVs significantly impacts the electricity network and the resulting demand profiles, particularly during peak periods. The surge in charging EVs during evenings can lead to grid instability and require additional investments in generation and transmission infrastructure \cite{leveque2007investments}.

The idea of community batteries promises to be a sound solution for increasing the integration of renewable energy in the power grid. These storage units, mainly located in the vicinity of ``communities,'' can store surplus renewable energy generated during peak production times and pass it back to the grid when demand is high. This helps to mitigate the variability and intermittency of renewable sources like solar and wind, ensuring more reliable and stable power system operations \cite{kalkbrenner2019residential}. 

The relevance of storage facilities for dealing with the intermittent fluctuation of PV resources is emphasised in several research publications. For instance, Fazlhashemi et al \cite{fazlhashemi2020day} describes an energy management model for distribution levels that are connected with energy storage while considering voltage stability, operational cost, and emissions as objectives of their models. Besides, the potential of batteries on power quality-index improvement was investigated in \cite{gholami2020energy}. On the same track, Sheidaei et al. \cite{sheidaei2021stochastic} develop a scenario-based scheduling strategy for energy storage in the presence of demand response programming. Congestion management through scheduling batteries in power systems was investigated by Yan et al \cite{yan2019robust}. Aryanezhad et al  \cite{aryanezhad2018management} focus on developing a scheduling approach for batteries using genetic algorithms to minimise voltage deviation and power loss.

Although these studies investigate certain crucial aspects of combining community batteries and renewable energy, their findings are at the utility-scale battery level that imposes massive expenditure on network administrators. In addition, it is debatable to increase the number of batteries in distribution networks as service providers may face challenges in managing their side effects as well as an excessive investment. It is important to note that the aforementioned publications mainly concentrated on the operational side of integrating batteries in power systems without rigorously analysing their optimal size and type for long-term planning.

Another significant aspect of utilising storage units in power systems can be community batteries, which are a promising remedy to capture the uncertainty associated with PV resources as they can be installed in the vicinity of rooftop solar panels. Some investigators have concentrated on the effects of community batteries in dealing with intermittency and variability of renewable generation. For example, Elkazaz et al.  \cite{elkazaz2021techno} focus on sizing community batteries with the aim of ancillary services to the market as well as managing the energy bills of customers. Similarly, Dinh et al. \cite{dinh2022optimal} develop an optimisation approach to size the community battery to participate in the local market. Although these papers cover interesting aspects of combining renewable resources and storage technologies, they have not considered EVs as part of the big picture. 

The optimal sizing of a community battery was the subject of several studies. Secchi et al. \cite{secchi2021multi} develop a multi-objective framework based on the non-dominated sorting genetic algorithm (NSGA-II) to optimally size community batteries. This investigation is beneficial as the community batteries were considered to manage the voltage fluctuations of distribution networks, but it suffers from using heuristic algorithms that cannot guarantee global solutions and expensive computational burden. Along similar lines, the effectiveness of community batteries in voltage management of distribution networks was discussed in \cite{alrashidi2022community}. This paper also investigates the size of community batteries to support the upstream network while it does not consider EV demands in the model.  The presence of EVs offers an opportunity for a more flexible power system. For example, smart charging technologies can be implemented to optimise charging schedules based on supply and demand profiles and individual preferences. By charging EVs during off-peak hours and using renewable energy resources, they can help reduce peak demand and contribute to grid stability.

Considering the rapid growth of EV utilisation, a plausible extension of the literature reviewed so far is the combination of renewable resources, in particular solar generation in this study, with community batteries that feed EV charging stations. However, several important components in this bundle require further investigation. For example, optimal long-term investment planning and sizing of the community batteries is the first step. This optimal planning and decision-making couple the historical renewable generation data with the charging and discharging patterns of community batteries and EVs. Accordingly, this study aims to develop a mixed-integer linear programming (MILP) model for effective sizing and scheduling of community batteries to capture the surplus generation of PV resources in a given community to supply the demand for EVs. Besides, it is necessary to assess the optimal investment plan under uncertainties originating from both generation and demand sides by a robust technique, without increasing the complexity of the MILP model.

Information gap decision theory (IGDT) establishes the foundation for rational decision-making in severe uncertainty. It proposes deterministic models to represent uncertain situations without requiring a huge amount of information. This approach differs from probability-based approaches to describe an uncertain situation when an extensive amount of data or Bayesian probabilistic models can rely on expert opinion. For a complicated stochastic system with a considerably small amount of data without expert opinion, IGDT serves as a suitable tool for decision-making \cite{colyvan2008probability}. IGDT investigates decisions (solutions) that are robust considering uncertainty. The robustness of a solution is quantified by the difference between the reward associated with the uncertain situation and a user-defined threshold \cite{ben2006info}.

Several studies have been carried out to tackle the problem of decision-making in the energy field in the presence of uncertainty. These approaches include methods such as fuzzy decision-making \cite{liang2011volt} for a power distribution system, and stochastic programming \cite{gholami2022multi} for EV charging stations. However, these methods have some limitations. For example, stochastic programming used in \cite{gholami2022multi} needs probability distributions of unknown parameters which are computationally expensive. IGDT for energy applications has been presented in \cite{gholami2023risk} and \cite{gholami2022risk} to address these obstacles. The IGDT consists of two main functions called the robustness and opportunity functions. In our investment and sizing problem, while the robustness function tries to increase the tolerance against risk by increasing the value of investments, the opportunity function tends to limit the exposure to uncertainties, which frequently leads to lower investments. One advantage of this approach is that it gives options to decision-makers based on their preferred level of risk from the two competing robustness and opportunity viewpoints.

In the context of renewable energy systems, where variables such as PV generation and EV demand exhibit considerable uncertainty, IGDT offers a compelling approach. The unpredictable nature of PV resources, coupled with the dynamic patterns of EV charging, necessitates decision-making strategies that can thrive in the face of evolving and unpredictable conditions. Wang et al. \cite{tostado2023information} and Li et al. \cite{correa2019comparative} stand as pioneers in integrating IGDT into scenario-based robust optimisation methods, emphasising its potential to strengthen optimisation models against unforeseen variations in renewable energy generation and EV usage. This integration introduces a decision-making framework explicitly designed to accommodate uncertainties, enhancing the adaptability of community battery systems. The incorporation of scenario-based approaches combines well with the unpredictable nature of PV generation and the variable demands imposed by EVs.

By integrating IGDT principles into our investment planning framework, we empower decision-makers with the necessary tools to navigate the complex landscape of renewable energy integration and electric vehicle (EV) charging along with community batteries. IGDT becomes not just a theoretical underpinning but a practical cornerstone in shaping resilient and future-ready strategies for community battery utilisation. The unique perspective offered by the main contribution of this work, opens avenues for addressing complexities associated with incomplete information, steering the optimisation of community battery deployment toward a sustainable and adaptive future. We anticipate continued advancements in renewable energy technologies and the increasing adoption of EVs, hence, this research highlights IGDT as a pivotal element in strengthening decision-making processes and ensuring the efficacy of community battery systems in dynamic and unpredictable environments.

This paper investigates the feasibility and benefits of using community batteries to collect excess power from solar panels and charge EVs. Establishing a sustainable and stable power source for charging EVs by aggregating the energy generated by solar panels within a community stored in batteries alleviates the challenges of integrating renewable energy into the electricity grid and reduces peak-period electricity demand. This method has several advantages, including load balancing, peak shaving, and improved grid stability. In addition, the intermittent nature of PV generation, uncertain electricity demand, and uncertain charging/discharging patterns of community batteries and EVs are captured using IGDT to address specific challenges encountered in community battery deployment, including real-time energy management, grid interactions, and the integration of emerging technologies.

The main contributions of this research can be outlined as follows:
\begin{itemize}
    \item Investigating the economic feasibility of community battery sizing by considering capital investment and operational costs. In more detail, a sophisticated MILP framework is developed for long-term strategic planning of community battery deployment inside a photovoltaic (PV)-rich region. The major goal of this framework is to efficiently fulfill the area's rising demand for electric vehicles (EVs) by identifying the optimal batter type, size, and investment plan. 
    \item A second critical consideration of our approach is the stability of the upstream grids. To achieve this, a power balancing constraint, which is especially important in circumstances where no surplus power is injected into upstream grids, is incorporated to restrict uncoordinated power injection of PV source to the grid. In other words, power injection to the grid is limited so that the unpredictable generation of PV resources does not impose any instability on the grid.
    \item To deal with uncertainties in the data from solar generation and EV charging demand, the IGDT framework is adapted to ensure a robust solution for investment planning and sizing. The IGDT is a nonparametric technique that does not need any information about the distribution of uncertain parameters.
\end{itemize}

The rest of this paper is structured as follows. Section \ref{sec_II} describes the MILP optimisation model for combining PV generation, community battery, and EVs. In Section \ref{sec_III}, the numerical results for the deterministic version of the problem are discussed. Sections \ref{sec_IV} and \ref{sec_V} demonstrate the IGDT framework for dealing with uncertainty and risk-based solution methods. 
Finally, the results and the discussion are provided in Section \ref{sec_Results}, which is then followed by the conclusion in Section \ref{sec_Conclusion}.

\section{The Mathematical Model} \label{sec_II}
\begin{figure}[!h]
  \centering
\includegraphics[width=\linewidth]{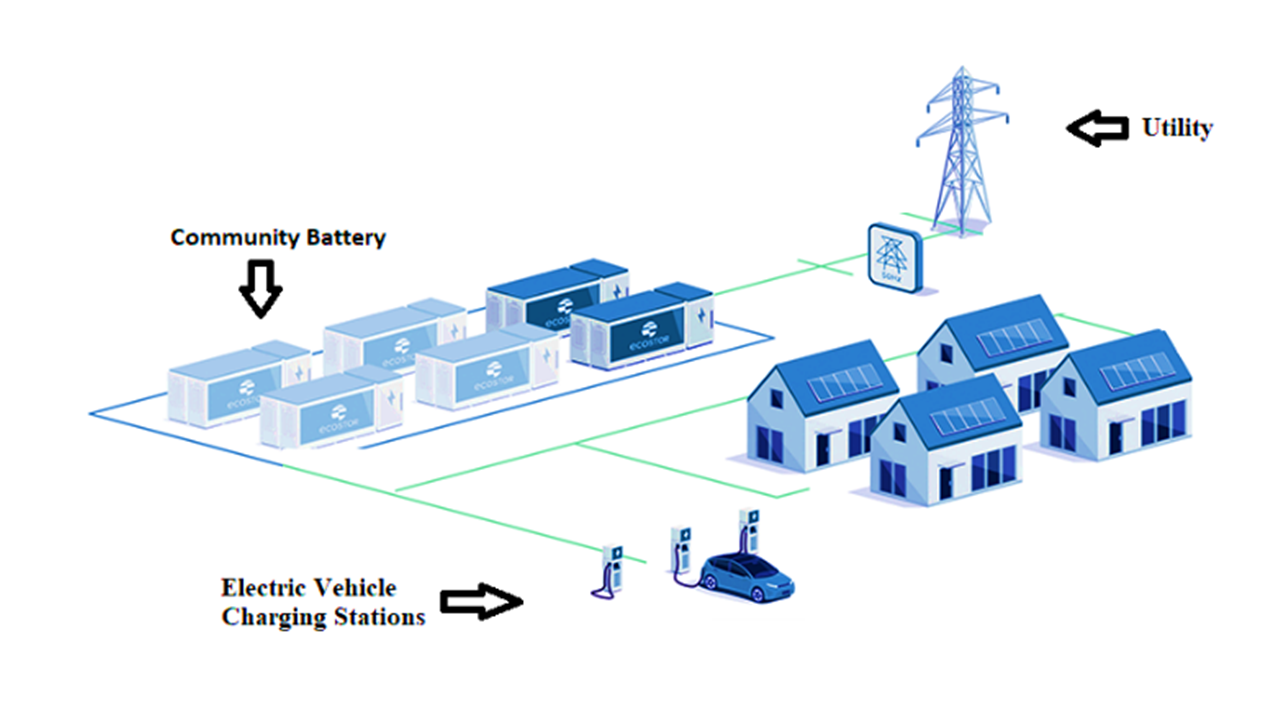}
  \caption{Integration of community battery in a residential area.}
  \label{comBat1_struc}
\end{figure}

A typical topology of a bundle of solar PV generation, community batteries, and an EV charging station in a sketch of a power system is depicted in Fig. \ref{comBat1_struc}. As mentioned earlier, the community battery accomplishes two main tasks: reducing the requirement for power exchange with the utility company by capturing excess power produced by PV resources and supplying stored power to charge electric cars. This combination not only improves the utility grid's stability by reducing the inherent uncertainties brought on by PV resources but also has financial benefits (discussion in Section~\ref{sec_Conclusion}). In this section, we develop a MILP model for finding an appropriate long-term investment plan and optimal sizing of community batteries for an area that includes electric vehicle demand, household demands, and rooftop PV generation.

The MILP model includes constraints on PV generation, battery charging/discharging behaviour, EV charging, and a power balance with an objective function focusing on minimising the capital investment and operation costs for a given planning horizon of $15$ years. The components of the MILP model including the indices, the parameters, and decision variables are shown in Table \ref{model_comp}.

\renewcommand{\arraystretch}{1.5}
\begin{table}
    \caption{The components of the MILP model}
    \centering
    \scalebox{0.6}{
    \begin{tabular}{|c|l|}
        \hline 
        \multicolumn{2}{|l|}{\textbf{Indices}} \\ \hline 
        $y \in \{1, \ldots, Y\}$ & year \\ \hline 
        $q \in \{1,2,3,4\}$ & quarter \\ \hline
        $d \in \{1, \ldots, T_q\}$ & days in a quarter $T_q \in \{90, 91,92,93\}$ \\  \hline 
        $t \in \{1, \ldots T=24\}$ & time of the day\\  \hline 
         $b \in \mathbb{B} = \{1,2,4,8\}$ & battery types \\  \hline 
         \multicolumn{2}{|l|}{\textbf{Parameters}} \\ \hline
        $\gamma$ & temperature derating \\  \hline 
        $\gamma_{_y}$ & discounting factor for year $y$ \\  \hline 
        $\eta_{_\text{PV}}$ & efficiency of PV generating unit \\  \hline
        $\eta_{_\text{ch}}$ & battery charging efficiency\\  \hline 
        $\eta_{_\text{dis}}$ & battery discharging efficiency \\  \hline 
        $M_{b}^{\text{max}}$ & a big number \\  \hline
        $I^C$ & solar insulation \\  \hline 
        $I^{\text{STC}}$ & insulation at standard test conditions \\  \hline 
        $T^{\text{CELL}}$ & cell temperature \\  \hline
        $T^{\text{AMB}}$ & ambient temperature \\  \hline 
        $\text{NOCT}$ & normal operating cell temperature \\  \hline 
         \multicolumn{2}{|l|}{\textbf{Decision variables}} \\ \hline
       $\text{Cap}_{b}^y$ & battery capacity of type $b$ installed in year $y$\\  \hline
       $C_{b}^y$ & investment cost of battery type $b$ installed in year $y$ outlined in Table  \ref{tableII}\\  \hline 
       $P_{\text{utility}}^{y,p,d,t}$ & amount of purchased power from the utility (grid) \\ \hline
       $C_{\text{utility}}^{y,p,d,t}$ & cost of buying power from the grid\\  \hline
       $S_{b}^{y,q,d,t}$ & state of charge (SoC) of a battery of type $b$ \\  \hline
       $P_{\text{b,ch}}^{y,q,d,t}$ & power injected for charging a battery\\  \hline
       $P_{\text{b,dis}}^{y,q,d,t}$ & power discharged form a battery\\  \hline
       $\text{Cap}_{b}^{\text{max},y}$ & the maximum battery capacity installed until year $y$ \\  \hline
       $B_{b}^{y,q,d,t}$ & binary variable dealing with simultaneous charging and discharging of batteries \\  \hline
       $P_{\text{PV}}^{y,q,d,t}$ &  power generated by PV units \\  \hline
       $P_{\text{load}}^{y,q,d,t}$ & local load demand \\  \hline
       $P_{\text{EV,ch}}^{y,q,d,t}$ & electric vehicle demand \\  \hline
    \end{tabular}}
    \label{model_comp}
\end{table}

\renewcommand{\arraystretch}{1}

The model's objective function attempts to minimise the summation of the investment and operational costs.
\begin{align}
   \text{OBJ}_{\text{\tiny{NPV}}} & =\text{Capex}+\text{Opex} \label{obj1} \\
   \text{Capex}&=\sum\limits_{y=1}^Y\sum\limits_{b \in \mathbb{B}} \gamma_{_y}  \text{Cap}_{b}^y  C_{b}^y \label{obj2}\\
   \text{Opex}&=\sum\limits_{y=1}^Y \sum\limits_{q=1}^4 \sum\limits_{d=1}^{T_q}\sum\limits_{t=1}^T \gamma_{_y}  P_{\text{utility}}^{y,p,d,t}  C_{\text{utility}}^{y,p,d,t}\label{obj3}
\end{align}
where 
$$\gamma_{_y}=\frac{1}{(1+r)^y}, \;\; y \in \{1, \ldots, Y\}$$
is the discount factor with $r$ representing the inflation rate.  In addition, Capex and Opex are the total capital investment and operational costs, $\text{Capacity}_{b}^y$ and $ C_{b}^y$ are the battery capacity and its investment cost at $y^{th}$ year for battery type $b \in \{1, 2, 4, 8\}$, and $P_{\text{utility}}^{y,p,d,t}$ and $C_{\text{utility}}^{y,p,d,t}$ are the amount of purchased power and its corresponding cost from the utility, respectively. The quantity $T_q \in \{90, 91, 92, 93\}$ is the number of days in different quarters in a year. The battery investment cost for different types of batteries throughout the planning horizon is provided in Table \ref{table1} which is extracted from CSIRO's GenCost technical report \cite{CSIRObat}. Essentially, Capex consists of multiplying the capacity of the battery and the associated cost for different types of batteries summed over the planning horizon with a discounting factor. It is important to note that we only consider a single battery type for the whole planning horizon, and a mixture of different battery types is not allowed due to some technical issues associated with battery controllers. The operational cost, Opex, is the summation of hourly power bought from the grid times the purchase price summed over the planning horizon with the same discounting factor.

Furthermore, the discounting factor $\gamma_y$ in \ref{obj3} enforces the net present value (NPV) of the investment and operational costs. As the planning horizon covers 15 years in this study, with a positive inflation rate, the current value of the money at present is more than the value of the money in the future. NPV with a discounting factor is in charge of holding the change of money value in this formulation. For some applications of NPV in project scheduling and optimisation see \cite{thiruvady2019maximising,Thiruvady2014ps}.

There are several sets of constraints in this model each in charge of a specific aspect of combining renewable energy, community battery, and EVs. Constraint set (\ref{const1}) takes care of the state of charge (SoC) of the batteries.
\begin{equation*}
    S_{b}^{y,q,d,t}=  S_{b}^{y,q,d,t-1} +\eta_{\text{ch}}P_{\text{b,ch}}^{y,q,d,t} 
    -\frac{P_{\text{b,dis}}^{y,q,d,t}}{\eta_{\text{dis}}} 
\end{equation*}
\begin{equation}
    \forall y,q,d, t\ge 1, \label{const1}
\end{equation}
SOC represents the remaining amount of energy in a battery considering the maximum amount possible to store in that battery when it is fully charged. The SoC of a battery of type $b$ at time $y,q,d,t$ is related to the amounts of charging and discharging of the party and the SoC of the previous hour. The SoCs are enforced to remain between the boundaries associated with the capacity of the batteries as indicated by (\ref{const2}).
\begin{equation}
    0 \le S_{b}^{y,q,d,t} \le \text{Cap}_{b}^{\text{max},y} \hspace{2.8cm} \forall y,q,d, t, \label{const2}
\end{equation}
The initial SoC of a battery at the beginning of the day is set as the final SoC of the same battery at the end of the previous day which is handled by (\ref{const7}).
\begin{equation}
    S_{b}^{y,q,d,0}=S_{b}^{y,q,d,T} \hspace{3.5cm} \forall y,q,d. \label{const7}
\end{equation}

To ensure that both charging and discharging do not occur at the same time, we use the big M technique, applied to constraints (\ref{const3}) and (\ref{const4})
\begin{equation}
    0 \le P_{\text{b,ch}}^{y,q,d,t}  \le B_{b}^{y,q,d,t} M_{b}^{\text{max}} \hspace{2.3cm} \forall y,q,d, t, \label{const3}
\end{equation}
\begin{equation}
    0 \le P_{\text{b,dis}}^{y,q,d,t}  \le (1-B_{b}^{y,q,d,t}) M_{b}^{\text{max}} \hspace{1.2cm} \forall y,q,d, t, \label{const4}
\end{equation}

Constraints (\ref{const5}) and (\ref{const6}) restrict the amount of charging and discharging power. we consider four different types of battery based on CSIRO's GenCost technical report \cite{CSIRObat} i.e., 1, 2, 4, and 8 hours types. As $\text{Cap}_{b}^{\text{max},y}$ represents the maximum battery capacity installed until year $y$, these constraints limit the amount of power inflow and outflow from a battery type $b$ in each hourly time interval $y,q,d, t$ based on the battery type.
\begin{equation}
    0 \le P_{\text{b,ch}}^{y,q,d,t}  \le \frac{\text{Cap}_{b}^{\text{max},y}}{b} \hspace{2.8cm} \forall y,q,d, t, \label{const5}
\end{equation}
\begin{equation}
    0 \le P_{\text{b,dis}}^{y,q,d,t}  \le \frac{\text{Cap}_{b}^{\text{max},y}}{b} \hspace{2.8cm} \forall y,q,d, t, \label{const6}
\end{equation}
The power balance at each hourly time interval $y,q,d, t$ is modelled using constraint (\ref{const8}). These constraints make sure that the amount of power generated by the solar units, purchased from the grid, stored and discharged from the batteries, the household demand, and utilised for charging EVs follow the power balance consideration.
\begin{equation*}
    P_{\text{PV}}^{y,q,d,t}+P_{\text{utility}}^{y,q,d,t} -P_{\text{b,ch}}^{y,q,d,t}+P_{\text{b,dis}}^{y,q,d,t} -P_{\text{load}}^{y,q,d,t}-P_{\text{EV,ch}}^{y,q,d,t}=0 
\end{equation*}
\begin{equation}
     \forall y,q,d,t, \label{const8}
\end{equation}
A high injection of renewably generated electricity into the grid can impose instability problems in a power system, and hence it is crucial to limit the power traded between community batteries and the grid. To accomplish this, constraints (\ref{const9}) restrict batteries from selling power to the utility. In more detail, a positive value of $P_{\text{utility}}^{y,q,d,t}$ indicates that we are only allowed to import energy from the utility and do not export the power to the grid. However, if $P_{\text{utility}}^{y,q,d,t}$ takes negative values, we essentially allow exporting power to the grid. In other words, constraints (\ref{const9}) ensure that the battery cannot be discharged into the grid to decline uncontrollable power injection into it.

\begin{equation}
     P_{\text{utility}}^{y,q,d,t}\ge 0 \hspace{4cm} \forall y,q,d,t.\label{const9}
\end{equation}
The restriction of injecting power into the grid has several advantages, including managing uncertainties associated with renewable resources which imply that the upstream network does not encounter instability issues because of power trading with distributed local resources.

The PV generation component in constraint (\ref{const8}) is modelled by equations (\ref{eq16}) based on the historical temperature and solar irradiation data in specific regions \cite{combe2019cost}:
\begin{align}
    & P^{y,q,d,t}_{\text{\tiny{PV}}}=\eta_{\text{\tiny{PV}}}P^{y,q,d,t}_{\text{\tiny{R}},\text{\tiny{PV}}}\left (\frac{I^C}{I^{STC}} \right ) (1-\gamma(T^{CELL}-T^{STC})), \notag \\
    &  \forall y, q, d, t, \label{eq16} \\
    &T^{CELL}=T^{AMB}+\left (\frac{\text{NOCT}-20^\circ}{0.8W/m^2} \right )I^C, \notag
\end{align}
where, as indicated in Table \ref{model_comp}, $\eta_{\text{\tiny{PV}}}$ is the efficiency of PV resources, $P^{y,q,t}_{\text{\tiny{R}},\text{\tiny{PV}}}$ is the rating of PV units, and $I^{STC}$ and $I^C$ are for insolation at standard test conditions and solar insolation, respectively. In addition, $\gamma$ is temperature derating, $T^{STC}$ represents the cell temperature, and $T^{AMB}(^\circ C)$ is ambient temperature. In this model,  $\text{NOCT}$ refers to normal operating cell temperature. It determines the operational conditions of solar cells, which include the highest temperature reached by open-circuited cells in a module. Generally, it is set as $800 W/m^2$ irradiance, $20^\circ C$ temperature, and wind speed of $1 m/s$ \cite{sun2020new}.

\section{Deterministic Solution Method} \label{sec_III}
For the deterministic version of the problem of finding an optimal investment plan and battery type and size, an objective function that includes both investment and operational costs for 15 years is considered in (\ref{MINLP1}). The determination of the type and capacity of the community battery happens in a way that the investment and operational costs are minimised. The investment cost is the amount of money invested to buy the units of batteries, which is defined as the purchase cost times the capacity. The operational costs are calculated based on the expenditure associated with purchasing power from the utility. Consequently, the problem with an embedded NPV in the objective function to find an optimal battery installation planning can be defined as the following optimisation problem:
\begin{align}
    P_{\text{MILP}}: &  \text{ Min}\;\; \text{OBJ}_{\text{\tiny{NPV}}}  \label{MINLP1} \\
    & \text{Subject to}   \notag\\
    & \hspace{15mm} \text{Constraints    } \text{(\ref{const1}) to (\ref{eq16})}  \notag
\end{align}

After solving the problem (\ref{MINLP1}), the payback period is determined using the following equation to identify the specific year at which the break-even point is reached.
\begin{equation}
    \text{Profit}^y=\sum\limits_{k=1}^y \text{Opex}^k_{\text{NoBattery}}-\text{Opex}^{y}_{\text{WithBattery}}, \label{eq18}
\end{equation}
where
$\text{Opex}^k_{\text{NoBattery}}$ and $\text{Opex}^y_{\text{WithBattery}}$ represent the operational costs without and with the community battery, respectively.

\begin{figure}[h]
  \centering
\includegraphics[width=0.7\linewidth]{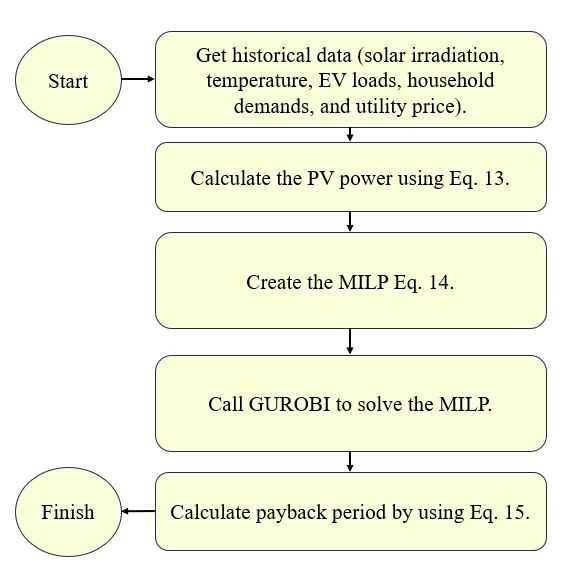}
  \caption{Flowchart of the solution approach to the problem (\ref{MINLP1}).}
  \label{Flowchart_struc}
\end{figure}

Fig. \ref{Flowchart_struc} demonstrates the solution approach to the problem (\ref{MINLP1}). After computing the PV generation using the historical data, a MILP model is created and solved using GUROBI\footnote{www.gurobi.com} solver, and the solution is fed into Eq. (\ref{eq18}) for the payback period calculations.

\section{Robust Solution Method: Information Gap Decision Theory} \label{sec_IV}
In the previous section, the deterministic formulation of the problem and its solution approach were discussed. However, it is impossible to find a promising investment plan without accounting for uncertainties in the problem including renewable energy generation and residential and EV demands. IGDT operates on the fundamental notion that decisions should exhibit robustness even when confronted with incomplete or imprecise information.

Theoretically, IGDT makes decisions by prioritising alternatives under uncertain factors, where there is a considerable gap between existing data and the complete information about uncertain components of a problem. It has two competing decision concepts called robustness and opportuneness. The robustness maximises the immunity against abrupt changes in uncertain parameters. The trade-off between the level of immunity and the quality of decisions is demonstrated via the robustness function. Opportuneness deals with better-than-expected decisions. In other words, robustness is the highest level for an uncertain parameter that warrants the required decision, and opportuneness is the lowest level for an uncertain parameter that enables an aspired decision \cite{ben2006info}.

\subsection{Robustness function}
IGDT uses several non-probabilistic uncertainty models to represent stochastic situations including fractional error and envelope-bound models. The envelope-bound model is used to handle both robustness and opportunity functions \cite{dai2019igdt}, \cite{samimi2019robust}.  Assume that the estimated and actual values of uncertainty parameters are shown by $\omega$ and $\Tilde{\omega}$ respectively. The aim is to construct the robustness function based on a reward function $R_{\text{OBJ}_{\text{NPV}}}$. It is defined as the maximum amount of uncertainty such that the minimum reward is greater than a predefined value. To attain this, let us consider $\alpha_r$ to be the uncertainty radius of uncertain parameters. The robustness function tends to maximise $\alpha_r$ by multiplying the uncertainty deviation factor ($\beta$) by the objective function defined in the deterministic case ($\text{OBJ}_{\text{NPV}}$).
\begin{align}
    & \text{Robustness function}= \notag \\
    &\max  \biggl\{\alpha_r: \max R_{\text{OBJ}_{\text{NPV}}} (\Xi,\Tilde{\omega})\le(1+\beta) \text{OBJ}_{\text{NPV}} (\Xi, \omega), \notag \\
    &\Tilde{\omega} \le (1\pm\alpha_r)\omega \biggr\}. \label{robustness}
\end{align}

Where, $\Xi$ are the decision variables (e.g. battery type, size, etc.). The expression $(1\pm\alpha_r)$  includes two operators ($\pm$) for which the minus sign and the plus sign are associated with PV generation and  EV demands, respectively, for reaching the worst-case scenario from the robustness function perspective.

The definition of the robustness function in (\ref{robustness}) constitutes a bi-level optimisation problem. However, it is feasible to streamline this problem into a single-level problem $P_{\text{rob}}$ based on the following reformulation.

\begin{align}
    P_{\text{rob}}: \text{ Max}\;\; \alpha_r & \label{Pprob} \\
    \text{Subject to}& \notag\\
    &  R_{\text{OBJ}_{\text{NPV}}} (\Xi,\Tilde{\omega})\le(1+\beta) \text{OBJ}_{\text{NPV}} (\Xi, \omega) \notag\\
    & \Tilde{\omega} \le (1\pm\alpha_r)\omega \notag\\
    &0 \le \alpha_r \le 1 \notag
\end{align}

\subsection{Opportunity function}
This function aims to find the best-case scenario because of uncertain parameter variations, such as an increase in PV generation and decreasing EV demand. Essentially, the opportunity of a decision (solution) is quantified as the minimum amount of uncertainty necessary for enabling the possibility of outcomes that exceed a critical value related to the reward function. Considering $\alpha_O$ as the uncertainty radius of uncertain parameters, the opportunity function is defined as follows.
\begin{align}
    & \text{Opportunity function}= \notag \\
    &\min \biggl\{\alpha_{\tiny{o}}: \min O_{\text{OBJ}_{\text{NPV}}} (\Xi,\Tilde{\omega})\ge (1-\beta) \text{OBJ}_{\text{NPV}} 
    (\Xi,\omega), \notag \\
    &\Tilde{\omega} \le (1\pm\alpha_o)\omega \biggr\} \label{opportunity}
\end{align}

From the opportunity function perspective, the expression $(1\pm\alpha_o)$ considers the minus sign and the plus sign for EV demands and PV generation respectively.

Again, similar to the definition of the robustness function, the definition of the opportunity function in (\ref{opportunity}) is based on a bi-level optimisation problem. The next expression shows the reformulation of it as a single-level optimisation problem $P_{\text{opp}}$. 
\begin{align}
     P_{\text{opp}}: \text{ Min}\;\; \alpha_o & \label{Popp} \\
    \text{Subject to}& \notag \\
    &  O_{\text{OBJ}_{\text{NPV}}} (\Xi,\Tilde{\omega})\ge(1-\beta) \text{OBJ}_{\text{NPV}}(\Xi, {\omega}), \notag \\
    & \Tilde{\omega} \le (1\pm\alpha_o)\omega \notag \\
    &0 \le \alpha_o \le 1 \notag
\end{align}

\section{Risk-based Solution Method} \label{sec_V}
In this section, the IGDT methodology is used to find robust solutions to the stochastic version of our problem. 
\subsection{Robustness Viewpoint}
 In our problem, the stochastic generation of solar PVs and EV demands are considered as the two uncertain parameters. Let $\alpha_r^{\text{PV}}$ and $\alpha_r^{\text{EV}}$ be the uncertainty radius of PV and EV, respectively. The extension of $P_{\text{rob}}$ for two uncertain radii is as follows:
\begin{align}
    P_{_R}: \text{ Max} & (\alpha_r^{\text{PV}}, \alpha_r^{\text{EV}}) & \label{Prob2} \\
    \text{Subject to}&  \notag \\
    &  R_{\text{OBJ}_{\text{NPV}}} (\Xi,\Tilde{P}_{\text{PV}}^{y,q,d,t}, \Tilde{P}_{\text{ev,ch}}^{y,q,d,t})\le  \notag \\
    &(1+\beta)\text{OBJ}_{\text{NPV}}(P_{\text{PV}}^{y,q,d,t}, P_{\text{EV,ch}}^{y,q,d,t}),  \;\;\forall y,q,d,t  \notag \\
    & \Tilde{P}_{\text{PV}}^{y,q,d,t} \le (1-\alpha_r^{PV})P_{PV}^{y,q,d,t}, \;\;\forall y,q,d,t \notag \\
    & \Tilde{P}_{\text{EV,ch}}^{y,q,d,t} \le (1+\alpha_r^{EV})P_{EV,ch}^{y,q,d,t}, \;\;\forall y,q,d,t  \notag\\
    & 0 \le \alpha_r^{PV} \le 1  \notag \\
    &0 \le \alpha_r^{EV} \le 1 \notag  \\
    & Constraints (4-12)  \notag
 \end{align}

\subsection{Opportunistic Viewpoint}
\begin{figure}[!h]
    \centering
    \includegraphics[width=\linewidth]{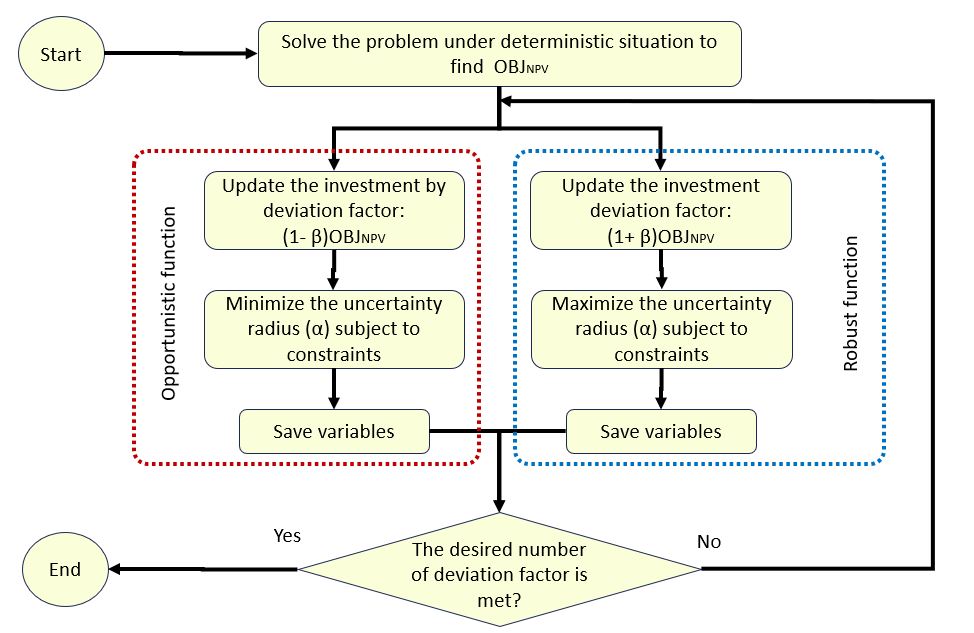}
    \caption{IGDT flowchart for handling uncertainty}
    \label{fig:IGDT}
\end{figure}

In contrast to the robustness function, with the opportunity function, we aim to find the minimum uncertainty radius to minimise the investment cost. This situation can be equivalently formulated as follows.

\begin{align}
    P_{_O}: \text{ Min} & (\alpha_o^{\text{PV}}, \alpha_o^{\text{EV}}) & \label{Popp2} \\
    \text{Subject to}&  \notag \\
    &  O_{\text{OBJ}_{\text{NPV}}} (\Xi,\Tilde{P}_{\text{PV}}^{y,q,d,t}, \Tilde{P}_{\text{ev,ch}}^{y,q,d,t})\ge  \notag \\
    &(1-\beta)\text{OBJ}_{\text{NPV}}(P_{\text{PV}}^{y,q,d,t}, P_{\text{EV,ch}}^{y,q,d,t}), \;\;\forall y,q,d,t \notag \\
    & \Tilde{P}_{\text{PV}}^{y,q,d,t} \le (1+\alpha_o^{PV})P_{PV}^{y,q,d,t}, \;\; \forall y,q,d,t \notag \\
    & \Tilde{P}_{\text{EV,ch}}^{y,q,d,t} \le (1-\alpha_o^{EV})P_{ev,ch}^{y,q,t}, \;\;\forall y,q,d,t \notag \\
    & 0 \le \alpha_o^{PV} \le 1 \notag  \\
    &0 \le \alpha_o^{EV} \le 1 \notag  \\
    & Constraints (4-12) \notag 
\end{align}
The flowchart in Fig \ref{fig:IGDT} demonstrates the process of solving the stochastic version using the IGDT approach for finding robust solutions. The deterministic problem (\ref{MINLP1}) is solved first to obtain the optimal cost represented by $\text{OBJ}_{\text{NPV}}$. Then, considering a user-defined value of $\beta$, the robust (\ref{Prob2}) and opportunistic problems (\ref{Popp2}) are solved to optimality to obtain the optimal values of $\alpha_r^{\text{PV*}}, \alpha_r^{\text{EV*}}, \alpha_o^{\text{PV*}}, \text{ and } \alpha_o^{\text{EV*}}$. These values assist a decision maker in finding optimal uncertainty radii for a given value of $\beta$. Solving problems (\ref{Prob2}) and (\ref{Popp2}) for different values of $\beta$ provide the risk and opportunity profiles of making decisions regarding investment options. These problems are solved for a set of $\beta$ values to create robust and opportunistic functions for uncertain radii.

\section{Results and Discussion} \label{sec_Results}
\subsection{Historical data and system setup}
\begin{figure*}[]
\centering
    \begin{subfigure}[t]{0.5\textwidth}
    \centering
    \includegraphics[width=0.8\textwidth]{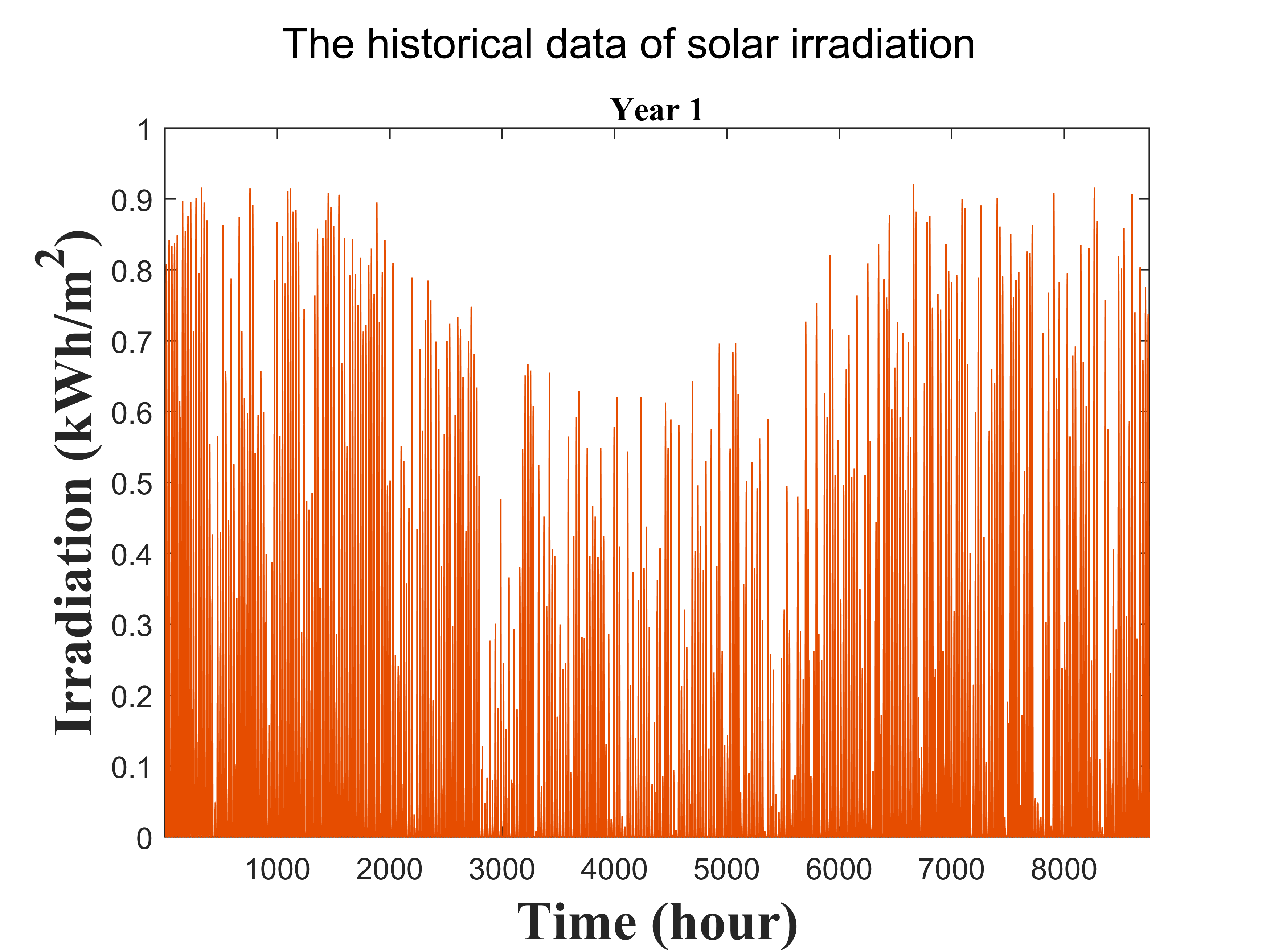}
    \caption{} \label{fig23}
\end{subfigure}\hfill
\begin{subfigure}[t]{0.5\textwidth}
    \centering
    \includegraphics[width=0.8\textwidth]{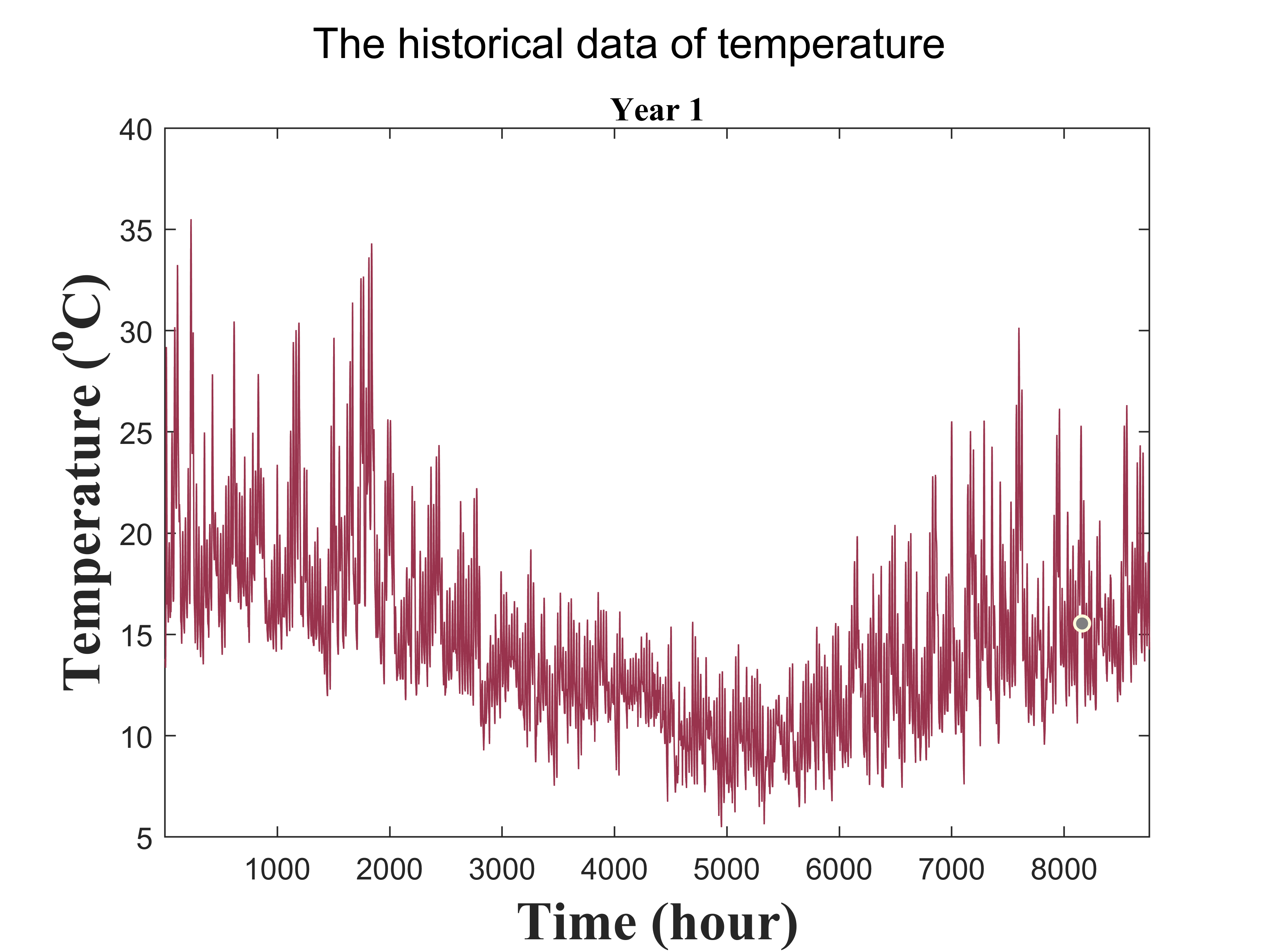}
    \caption{} \label{fig33}
\end{subfigure}\hfill
\begin{subfigure}[t]{0.5\textwidth}
    \centering
    \includegraphics[width=0.8\textwidth]{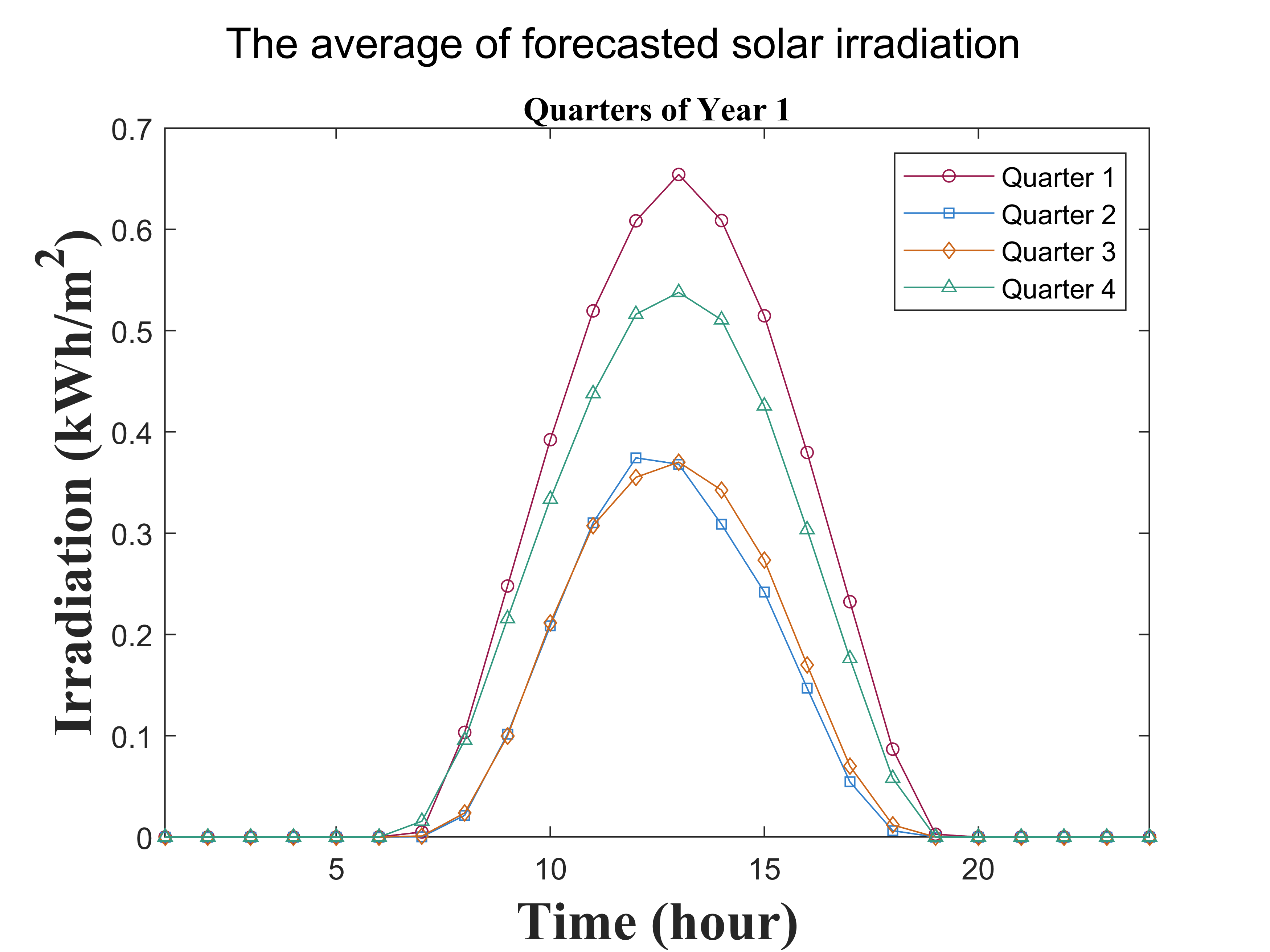}
    \caption{} \label{fig43}
\end{subfigure}\hfill
\begin{subfigure}[t]{0.5\textwidth}
    \centering
    \includegraphics[width=0.8\textwidth]{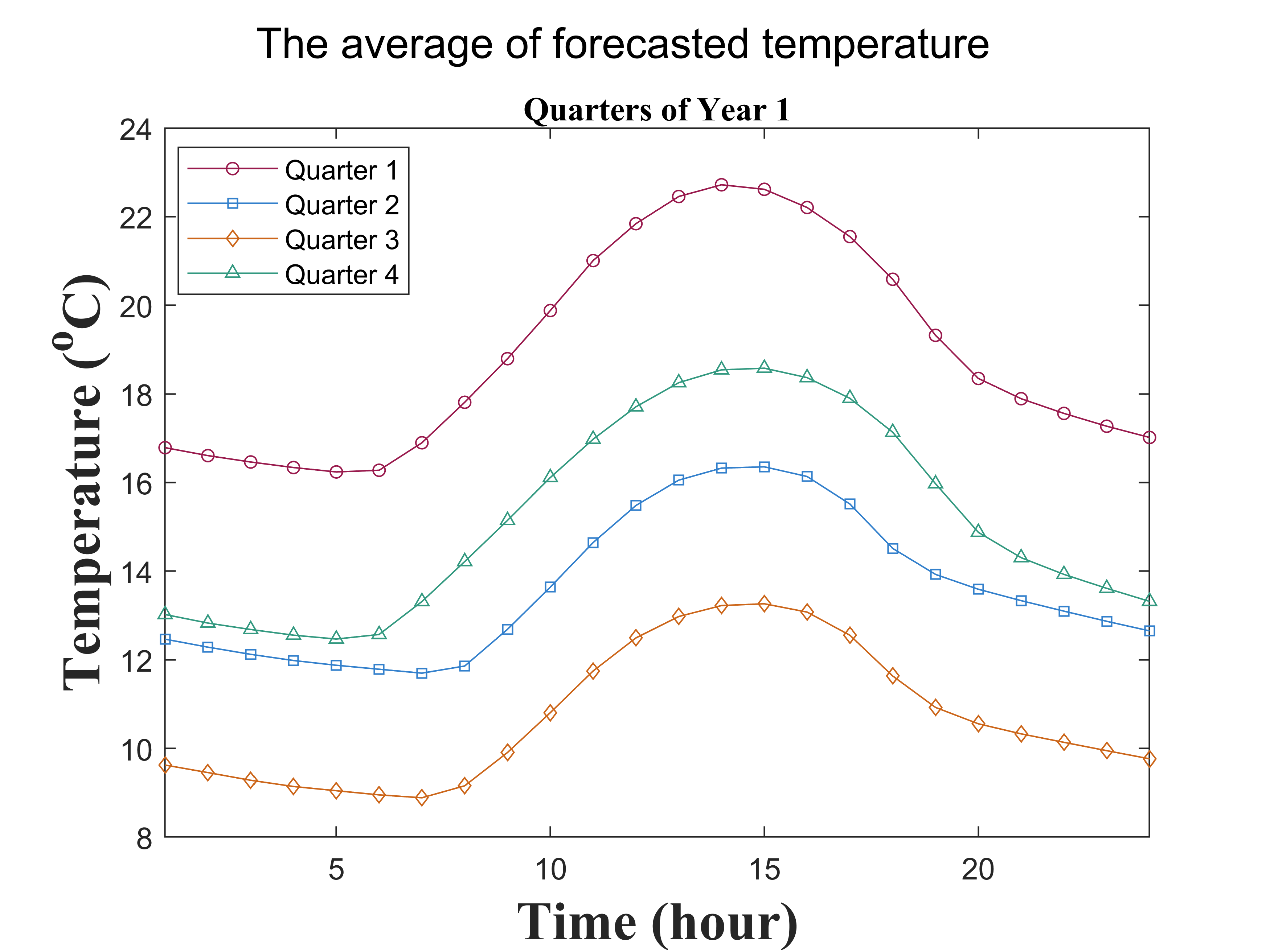}
    \caption{} \label{fig53}
\end{subfigure}\hfill
\begin{subfigure}[t]{0.5\textwidth}
    \centering
    \includegraphics[width=0.8\textwidth]{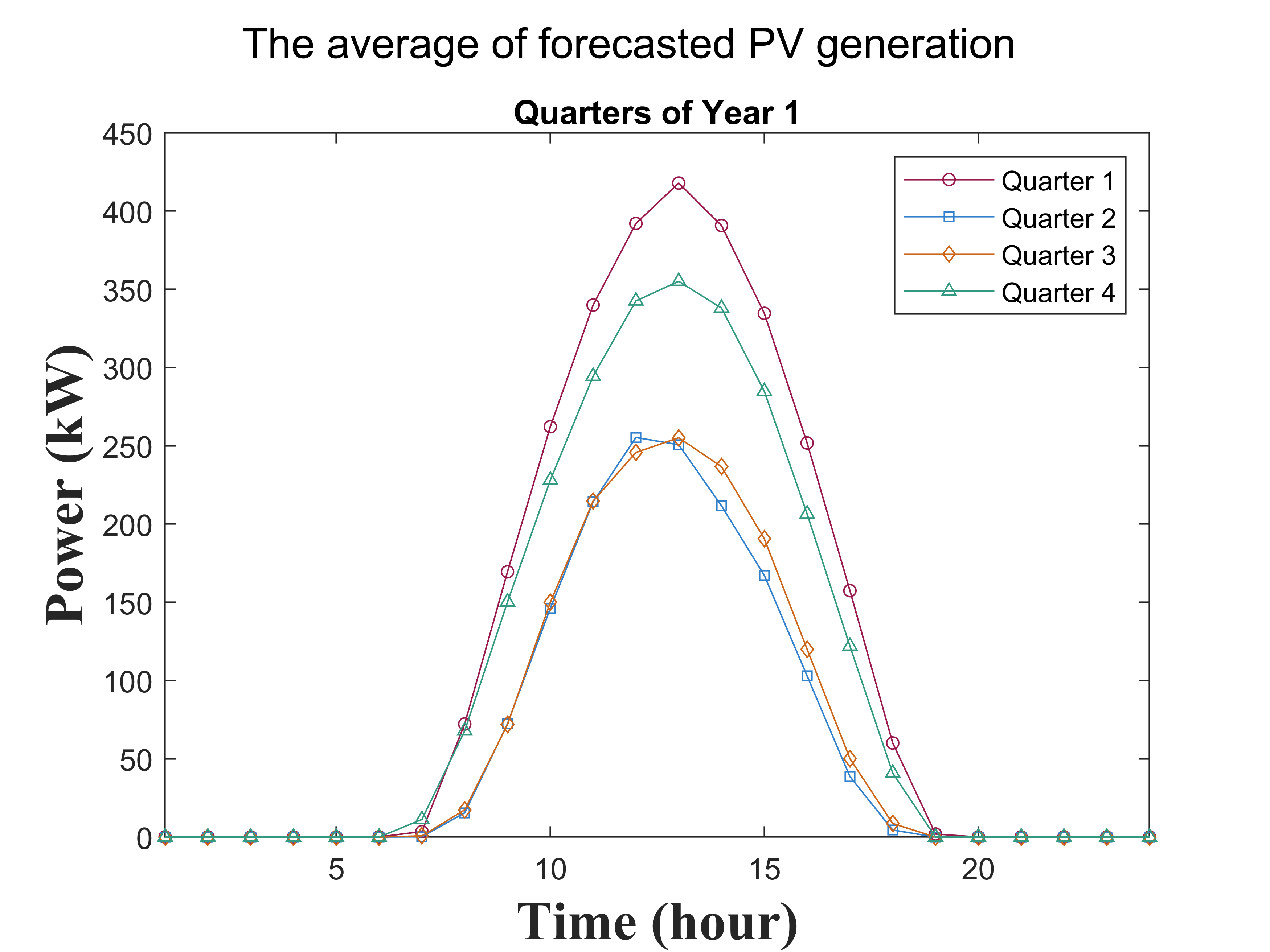}
    \caption{} \label{fig63}
\end{subfigure}\hfill
\begin{subfigure}[t]{0.5\textwidth}
    \centering
    \includegraphics[width=0.8\textwidth]{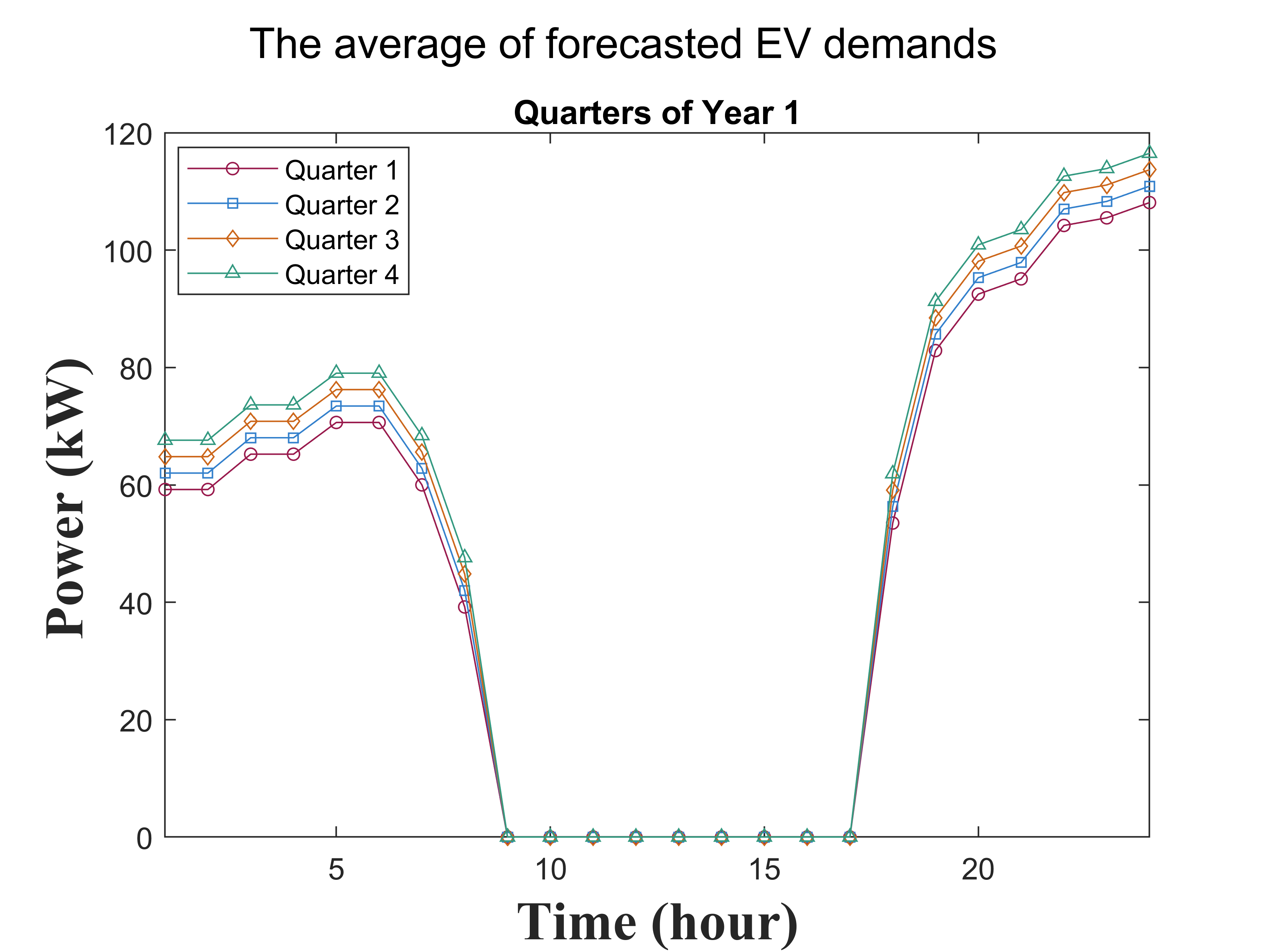}
    \caption{} \label{fig73}
\end{subfigure}\hfill
\begin{subfigure}[t]{0.5\textwidth}
    \centering
    \includegraphics[width=0.8\textwidth]{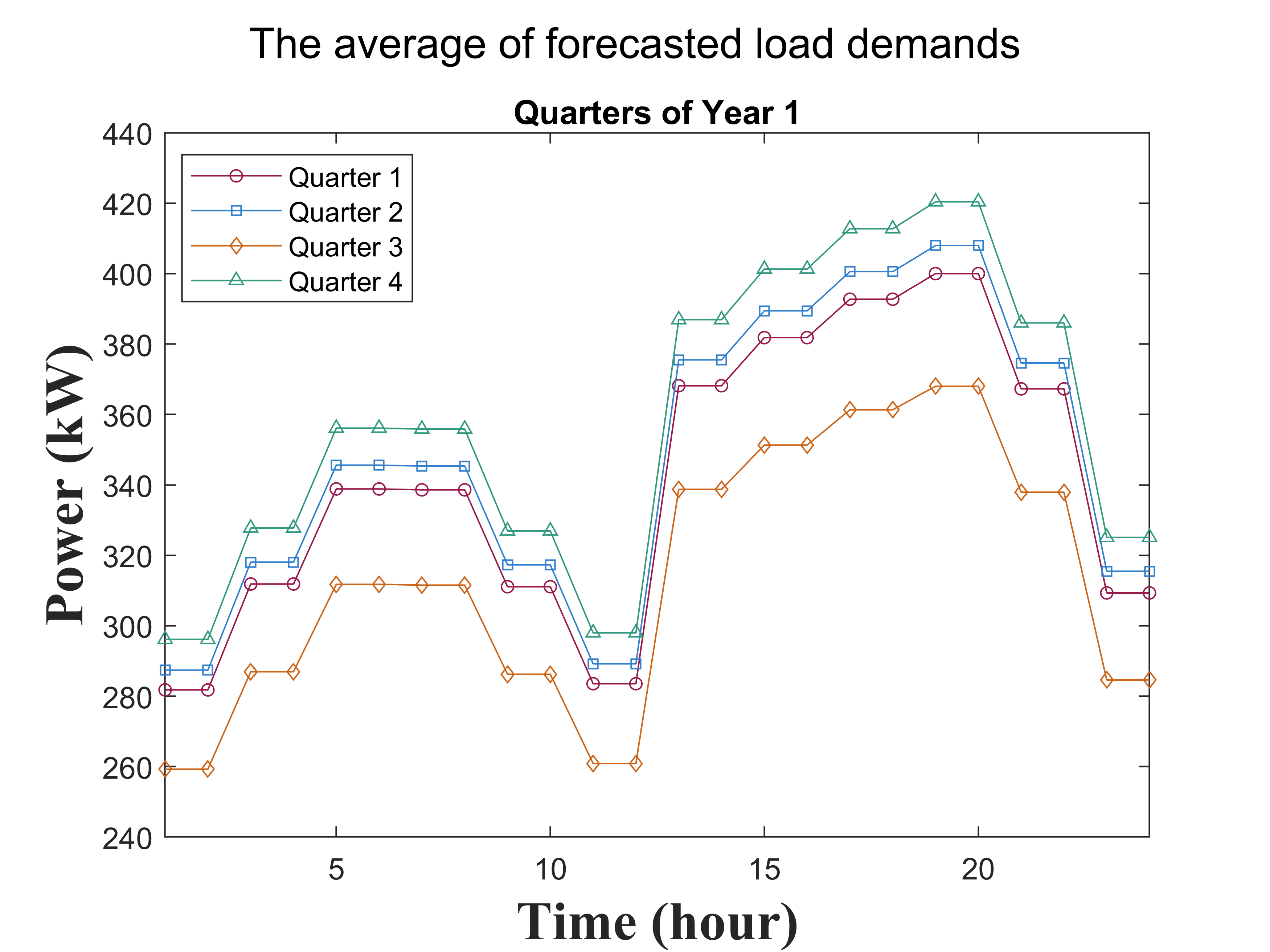}
    \caption{} \label{fig83}
\end{subfigure}\hfill
\begin{subfigure}[t]{0.5\textwidth}
    \centering
    \includegraphics[width=0.8\textwidth]{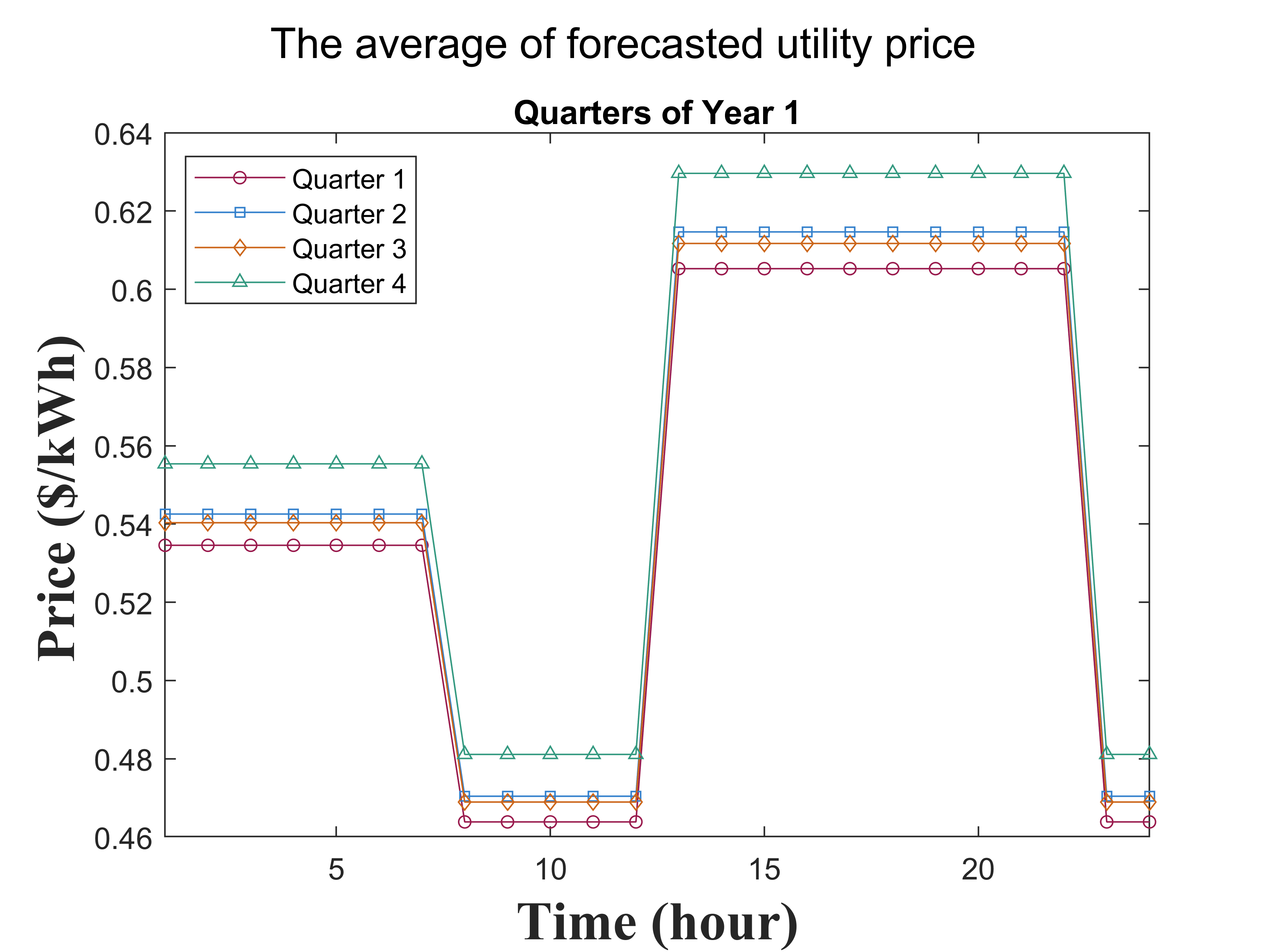}
    \caption{} \label{fig93}
\end{subfigure}\hfill
\caption{For the first year of the planning horizon, (a) historical data of solar irradiation; (b) historical data of temperatures;  (c) quarterly Average Solar Irradiation;  (d) quarterly average temperature;  (e) quarterly PV Resource Output; (f) Quarterly Average Electric Vehicle Demands;(g) quarterly average residential load demands;  (h) quarterly Utility Prices. The presentation of the data for the whole planning horizon is in Appendix \ref{appendix}.}
\label{fig:mylabel3}
\end{figure*}

To construct the optimisation models, a planning horizon of 15 years ($Y=15$) for community batteries is considered in this work. A long-term planning horizon, instead of a single year, is necessary to capture the growth in EVs and PV generation. The historical solar irradiation and temperature data of Geelong, Victoria, Australia are gathered and summarised in Figs. \ref{fig23} and  \ref{fig33}, respectively\footnote{Renewable energy, \url{https://www.renewables.ninja}}. The quarterly average of the irradiation and temperature data are depicted in Figs. \ref{fig43} and \ref{fig53}.  These average data are converted into solar-generated power using Eq. (\ref{eq16}) which are presented in Fig. \ref{fig63} and are used in the optimisation model for planning. Fig. \ref{fig73} shows the quarterly EV demand for 15 years, and the residential demand is depicted in Fig. \ref{fig83}. Fig. Also, \ref{fig93} shows the utility price for the planning horizon.  Note that Fig. \ref{fig:mylabel3} only shows the data for the first year of the planning horizon due to the readability of the plots. The plots of data incorporated in our models for the whole planning horizon are depicted in Appendix \ref{appendix}. 

The MILP optimisation problem is implemented in Python and solved using GUROBI \cite{gurobi}. A Personal i5 laptop with 8GB of memory is utilised to obtain the numerical results. Two case studies are considered to evaluate the effect of community batteries: 
\begin{itemize}
    \item \textbf{Case 1:} In this case, we do not include community batteries. Therefore, the EVs are charged using grid electricity and PV-generated electricity. 
    \item \textbf{Case 2:} In this case, community batteries are included.
\end{itemize}

\begin{table}[]
\large
\centering
\caption{The price of battery types (\$/kWh) based on CSIRO's GenCost technical report \cite{CSIRObat}.}
\label{table1}
\centering
\scalebox{0.7}{
\begin{tabular}{c|cccc}
\multicolumn{1}{l}{} & \multicolumn{4}{c}{\textbf{Battery Type}} \\ \cline{2-5}
\textbf{Year}         & \textbf{1h} & \textbf{2h} & \textbf{4h} & \textbf{8h} \\ \hline
2023 & 935         & 676         & 549         & 487         \\ \hline
2024 & 879         & 635         & 516         & 458         \\ \hline
2025 & 830         & 600         & 487         & 433         \\ \hline
2026 & 786         & 568         & 461         & 410         \\ \hline
2027 & 791         & 554         & 448         & 397   \\ \hline
2028 & 773         & 539         & 441         & 396         \\ \hline
2029 & 755         & 525         & 427         & 383         \\ \hline
2030 & 737         & 510         & 414         & 371         \\ \hline
2031 & 719         & 496         & 401         & 358         \\ \hline
2032 & 701         & 481         & 388         & 345         \\ \hline
2033 & 683         & 467         & 374         & 332         \\ \hline
2034 & 665         & 453         & 361         & 320         \\ \hline
2035 & 647         & 438         & 348         & 307         \\ \hline
2036 & 629         & 424         & 335         & 294         \\ \hline
2037 & 612         & 410         & 322         & 282              
\end{tabular}}
\end{table}

\subsection{Deterministic Evaluation}
An instance of the MILP problem (\ref{MINLP1}) is created by incorporating PV electricity generation values obtained from the historical data and the MILP is solved using GUROBI. It is important to note that only one battery type, either type 1, 2, 4, or 8, in the solution, and a mixture of different types is not allowed due to technical issues concerning battery controllers. The optimal solution to the problem is then translated into the best investment planning for a horizon of 15 years in this study. 
\begin{table*}[!h]
\caption{The battery capacity of different types built per year (kWh)}
\label{tableI}
\centering
\scalebox{0.65}{
\begin{tabular}{c|ccccccccccccccc}
 Bat Type              & Y1 & Y2 & Y3 & Y4 & Y5 & Y6 & Y7 & Y8 & Y9 & Y10 & Y11 & Y12 & Y13 & Y14 & Y15 \\ \hline
1 & 226   & 581   & 481   & 1728  & 0     & 362   & 183   & 198   & 158   & 777    & 0      & 0      & 0      & 0      & 0      \\  \hline
2 & 226   & 581   & 481   & 1667  & 79    & 344   & 206   & 175   & 158   & 1439   & 447    & 0      & 0      & 0      & 0      \\  \hline
4 & 434   & 544   & 1721  & 317   & 96    & 266   & 212   & 169   & 238   & 1359   & 515    & 219    & 0      & 403    & 275    \\  \hline
8 & 868   & 1087  & 782   & 820   & 927   & 879   & 653   & 431   & 0     & 2026   & 615    & 214    & 0      & 662    & 1466  
\end{tabular}}
\end{table*}

Tables \ref{tableI} -- \ref{tableIV} show costs and installed battery capacities corresponding to the optimal solution of the problem (\ref{MINLP1}) built using the given data. Table \ref{tableI} shows the yearly capacity installation of each battery type for 15 years. The level of capacities for each year is optimally decided to cope with increasing solar energy generation, EV charging demands, and household electricity load. Notably, all the necessary battery capacity for types 1 and 2 are installed in the first 10 and 11 years respectively, while for the others, it is gradually installed over the 15 year period. The total sum of installed capacities for types 1,2,4, and 8 are 3994, 5803, 6768, and 11430 kWh, respectively. The technical reason behind different installed capacity levels is that battery type 1, for example, has a more charging-discharging boundary ($P_{1,\text{ch}}, P_{1,\text{dis}} \le \frac{\text{Cap}_1}{1}$), as outlined in constrained (\ref{const5}) and (\ref{const6}), which implies that less battery capacity is needed. However, in other types, particularly in type 8, the boundary of battery charging-discharging is confined by ($P_{8,\text{ch}}, P_{8,\text{dis}} \le \frac{\text{Cap}_8}{8}$). As a result, more capacity is utilised. In addition, Table \ref{tableII} provides the monetary plan for yearly investment for each battery type. Considering the summation of each row for installed capacity in Table \ref{tableI} and Capex in Table \ref{tableII}, the minimum capacity installed is achieved by battery type 1 with the value 3,994 kWh, and the minimum investment occurred for battery type 4 with the value \$2,587,000.
\begin{table*}[!h]
\caption{The Capex per year (\$1000/year)}
\label{tableII}
\centering
\scalebox{0.65}{
\begin{tabular}{c|ccccccccccccccc}
 
Bat Type              & Y1 & Y2 & Y3 & Y4 & Y5 & Y6 & Y7 & Y8 & Y9 & Y10 & Y11 & Y12 & Y13 & Y14 & Y15 \\  \hline
1 & 210 & 497 & 379 & 1260  & 0        & 246 & 118 & 122 & 93 & 435 & 0        & 0        & 0      & 0        & 0        \\  \hline
2 & 152 & 359 & 274 & 879 & 39 & 163 & 93 & 74 & 64 & 554 & 162 & 0        & 0      & 0        & 0        \\  \hline
4 & 238 & 273 & 797 & 135 & 38 & 103 & 77 & 58 & 78 & 422   & 150 & 59 & 0      & 97 & 62 \\  \hline
8 & 422   & 485 & 321 & 311 & 333 & 307 & 215 & 134 & 0        & 559 & 159 & 52 & 0      & 141 & 292
\end{tabular}}
\end{table*}

Tables \ref{tableIII} and \ref{tableIV} provide detailed information regarding the operational cost, Opex, for the two scenarios of a power system with or without battery components, as outlined in Case 1 and Case 2. Unsurprisingly, the Opex for a power system with batteries is considerably lower than the case without batteries. The capability to store surplus solar generation and subsequently discharge it back to the grid during peak utility prices emerges as the driving force behind these cost reductions. In other words, when a community battery is incorporated, the reliance on purchased power from the utility diminishes, which leads to a reduction in Opex.
\begin{table*}[!h]
\caption{The Opex (\$1000/year) per year when battery is not considered.}
\label{tableIII}
\centering
\scalebox{0.62}{
\begin{tabular}{c|ccccccccccccccc}
Bat Type              & Y1 & Y2 & Y3 & Y4 & Y5 & Y6 & Y7 & Y8 & Y9 & Y10 & Y11 & Y12 & Y13 & Y14 & Y15 \\ \hline  
1 & 1496 & 1455 & 1504 & 1509 & 1559 & 1609 & 1900 & 1994 & 2101 & 2199 & 2291 & 2358 & 2478 & 2508 & 2620 \\  \hline
2 & 1496 & 1455 & 1504 & 1509 & 1559 & 1609 & 1900 & 1994 & 2101 & 2199 & 2291 & 2358 & 2478 & 2508 & 2620 \\  \hline

4 & 1496 & 1455 & 1504 & 1509 & 1559 & 1609 & 1900 & 1994 & 2101 & 2199 & 2291 & 2358 & 2478 & 2508 & 2620 \\  \hline
8 & 1496 & 1455 & 1504 & 1509 & 1559 & 1609 & 1900 & 1994 & 2101 & 2199 & 2291 & 2358 & 2478 & 2508 & 2620
\end{tabular}}
\end{table*}

The sum of the total battery investment cost and the operational cost for the planning horizon is summarised in Table \ref{tableV} for both scenarios. When there is no battery installed in the system, there is no investment cost, and the total operational cost throughout the planning horizon is \$29,587,344. As there is no investment, there is no profit as well. In the case of installing batteries in the system, depending on the battery type the total cost as the summation of Capex and Opex for 15 years, and the amount of profit are presented in the table. Out of four battery types, the optimal choice is the battery type 4, which incurs \$ 23,194,628 of total cost with a profit of \$ 6,392,716. The amount of yearly profit for each battery type is calculated using Eq. (\ref{eq18}) and summed up.
\begin{table*}[!h]
\caption{The Opex (\$1000/year) per year when battery is considered}
\label{tableIV}
\centering
\scalebox{0.62}{
\begin{tabular}{c|ccccccccccccccc}
Bat Type              & Y1 & Y2 & Y3 & Y4 & Y5 & Y6 & Y7 & Y8 & Y9 & Y10 & Y11 & Y12 & Y13 & Y14 & Y15 \\ \hline
1 & 1476 & 1361 & 1349 & 1220 & 1201 & 1129 & 1366 & 1385 & 1493 & 1462 & 1455 & 1472 & 1602 & 1577 & 1680 \\  \hline
2 & 1476 & 1361 & 1349 & 1221 & 1201 & 1129 & 1365 & 1385 & 1493 & 1428 & 1397 & 1402 & 1564 & 1460 & 1561 \\  \hline
4 & 1468 & 1355 & 1297 & 1220 & 1199 & 1129 & 1365 & 1385 & 1492 & 1428 & 1394 & 1387 & 1562 & 1406 & 1503 \\  \hline
8 & 1460 & 1336 & 1321 & 1235 & 1200 & 1125 & 1356 & 1375 & 1484 & 1419 & 1382 & 1375 & 1551 & 1392 & 1455
\end{tabular}}
\end{table*}

\begin{figure}[!h]
  \centering
\includegraphics[width=\linewidth]{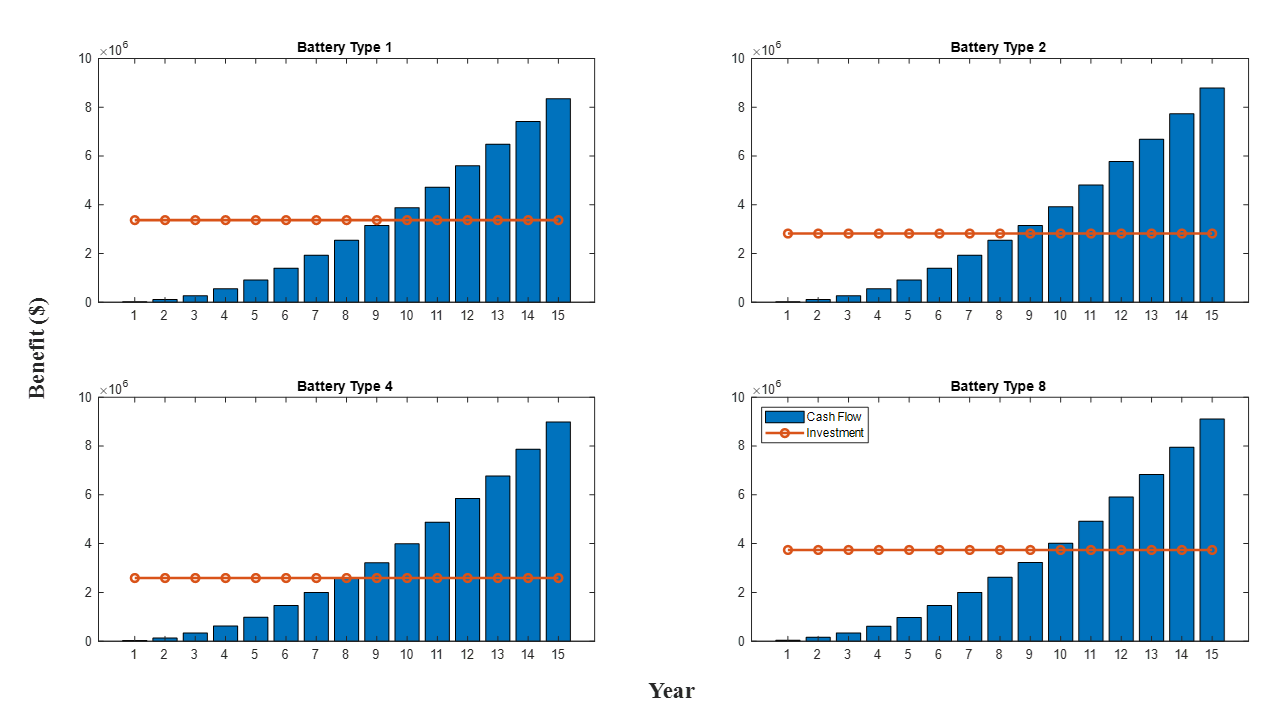}
  \caption{ Cashflow Analysis for Various Battery Types. }
  \label{comBat1}
\end{figure}

Cash flow is utilised to determine the number of years required for the investment to break even. This concept is also visually represented in Fig. \ref{comBat1}. The Capex value is shown by a red line, and the yearly Opex is depicted using blue bars. The area under the bars above the investment line shows the amount of profit that is the maximum for battery type 4. An investor begins to realise the profit after 8 years, and it takes another 8 years for them to retrieve their initial investment.
\begin{table}[h]
\caption{Total cost of 15 years for different battery types versus no battery}
\label{tableV}
\centering
\scalebox{0.7}{
\begin{tabular}{c|cc}
Method   & Total Cost (\$/15 years) & Profit (\$/15 year) \\ \hline
No battery     & \$29,587,344             & -                   \\
Battery type 1 & \$24,601,350             & \$4,985,994         \\
Battery type 2 & \$23,616,520             & \$5,970,824         \\
Battery type 4 & \$23,194,628             & \$6,392,716         \\
Battery type 8 & \$24,211,920             & \$5,375,424        
\end{tabular}}
\end{table}

\subsection{Risk-based Evaluation}
The deterministic solution of the problem (\ref{MINLP1}) provides us with an optimal solution. It helps us to find the most appropriate battery type, establish the investment planning for a horizon of 15 years in this study, and compute the operational costs. In this section, we focus on the implications of the IGDT framework on risk-based decision-making by constructing robustness and opportunity functions for different deviation factors. For a given deviation factor value $\beta$, the decision-maker is interested in obtaining insights about the amount of increase or decrease in the total cost and finding out the corresponding radii of uncertainty. As long as the uncertain parameters fluctuate in the radii of uncertainty, the solution will remain feasible. For constructing robustness and opportunity functions, after solving the deterministic problem and obtaining the optimal value $\text{OBJ}_{NPV}$, we solve the robustness problems (\ref{Prob2}) and opportunistic problem (\ref{Popp2}) for different values of $\beta \in \{\beta_0, \ldots, \beta_n\} \subset (0,1]$. These values play the role of discounting coefficients in a way that the optimal solutions for problems (\ref{Prob2}) and (\ref{Popp2}) are not worse than the optimal solution of the original problem this discounting coefficient. For example, by setting $\beta=0.1$, we are interested in finding the optimal values for uncertainty radii $\alpha_r^{\text{PV*}}, \alpha_r^{\text{EV*}}, \alpha_o^{\text{PV*}}, \text{ and } \alpha_o^{\text{EV*}}$ when a 10\% deviation is allowed form the optimal total cost $\text{OBJ}_{NPV}$. These optimal values for uncertainty radii are incorporated to construct robustness and opportunity functions. The number of different values for the deviation factor $\beta$ and its distribution is decided by the decision-maker.

Fig. \ref{comBat1} is a typical outcome of solving problems (\ref{Prob2}) and (\ref{Popp2}) for a set of $\beta$ values. The $x$-axis represents the total cost, and the $y$-axis is the radius of uncertainty. It shows the objective value of the optimal deterministic problem multiplied with different values of $\beta$ as $(1\pm\beta) \text{OBJ}_{\text{NPV}}$. The blue lines represent the robustness functions, and the red lines are the opportunistic functions for EV (the top plot) and PV (the bottom plot). For each plot, the values of $\alpha_r$ and $\alpha_o$ corresponding to a particular $\beta$ value are read using the blue and red lines correspondingly. When $\beta=0$, we deal with the deterministic case, and  $\alpha_r^{\text{EV}}=\alpha_o^{\text{EV}}=0$.  
 
\begin{figure}[!h]
    \centering
    \includegraphics[width=\linewidth]{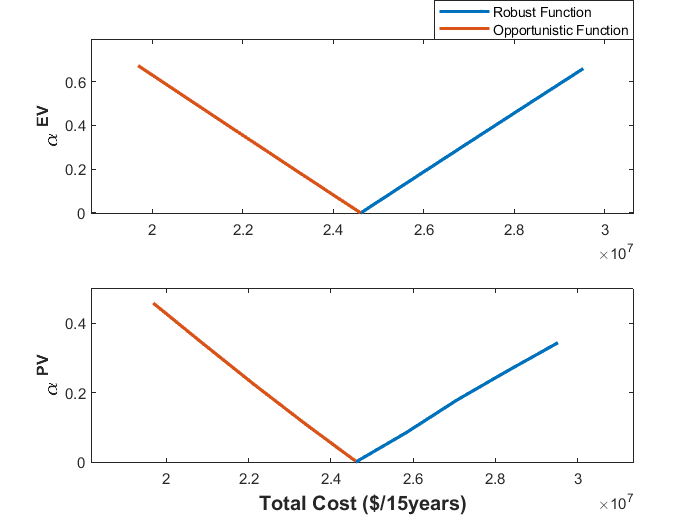}
    \caption{The uncertainty radiuses from both robust (in blue) and opportunistic (in red) viewpoints. For a given $\beta$ the $\alpha_r^{\text{EV}}$ and $\alpha_r^{\text{PV}}$ are understood using the blue lines, and the $\alpha_o^{\text{EV}}$ and $\alpha_o^{\text{PV}}$ are obtained using the red lines}
    \label{fig:IGDTplot}
\end{figure}

To examine robustness in this study, illustrated by Fig. \ref{comBat1}, investors may want to allocate additional resources to ensure a resilient system in response to varying uncertainty due to changes in solar renewable generation or EV load. This strategy strengthens the model against an increase in EV demand and a decrease in PV resources, prioritising worst-case scenarios.  In contrast, the opportunity function presents an opposite action choice. Here, investors seek to reduce investment when the generation of PV resources increases while EV demand decreases. This approach optimises scenarios where favourable conditions permit reduced investment. The robustness function is particularly relevant in practical applications where the worst-case scenario carries more weight (Fig. \ref{fig:IGDTplot}). This slight evaluation through IGDT provides a comprehensive risk-based perspective, enriching the decision-making process for community battery deployment. 

\begin{figure}[!h]
    \centering
    \includegraphics[width=\linewidth]{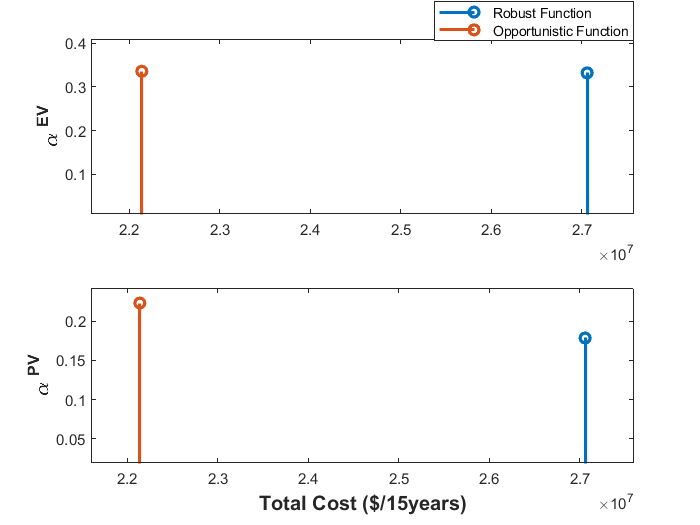}
    \caption{For a user-defined value of $\beta=0.1$, the corresponding uncertainty radius values are $\alpha_r^{\text{EV}}=0.3$ and $\alpha_r^{\text{PV}}=0.18$, $\alpha_o^{\text{EV}}=0.35$, and $\alpha_o^{\text{PV}}=0.25$, respectively.}
    \label{fig:IGDT10perc}
\end{figure}

Fig. \ref{fig:IGDT10perc} provides an example of computing the radii of uncertainty for a given value of $\beta=0.1$ in Fig. \ref{fig:IGDTplot}. Drawing a vertical line at $\beta=0.1$ which is equivalent to $1.1\text{OBJ}_{\text{NPV}}$, in both curves for EV and PV, gives us the radius of uncertainty of  $\alpha_r^{\text{EV}}=0.3$ and $\alpha_r^{\text{PV}}=0.18$ form the robustness point of view.

\section{CONCLUSION} \label{sec_Conclusion}
The increasing popularity of EVs originates from their remarkable capacity to significantly decrease greenhouse gas emissions and reduce dependence on fossil fuels. However, this paradigm shift from traditional fossil fuel-based vehicles towards EVs is not without challenges, particularly in the existing electrical power systems which are struggling with the technical complexities of high integration of scattered solar resources within communities. Our work addresses this conjunction between renewable energy and EVs by proposing a strategic planning approach aimed toward evaluating the long-term effectiveness of community battery sizing over a 15-year horizon.

Intending to optimise community battery deployment, we considered the economic implications by incorporating the real price of various battery types. Our findings underscore the efficacy of battery type 4, revealing it as the most cost-efficient option with a faster break-even trajectory compared to other types. Beyond the economic landscape, our research illustrates the profound environmental impact, representing that battery communities have the potential to generate substantial cost savings, amounting to \$6 million, in contrast to scenarios without their inclusion.

Furthermore, we recognise the inevitability of uncertainties in the dynamic interplay between EV demand and PV generation. To strengthen our planning strategy against uncertainties, we applied the Information-Gap Decision Theory technique, examining the model through the lens of both robust and opportunistic functions. Notably, the robust function demonstrated a tendency to incur additional costs in exchange for minimising vulnerability to deviations in EVs and PVs. Conversely, opportunistic functions adopted a cost-saving approach, highlighting the trade-offs inherent in optimising community battery deployment under uncertainties.

As we move towards future investigations, we plan a thorough exploration of the complicated realm of demand response programming. This endeavour will add a layer of sophistication to our model, offering a deeper understanding of the synergies between demand response strategies and the utilisation of community batteries. Through this exploration, we aim to enhance the resilience and adaptability of our planning approach, paving the way for a more comprehensive, dynamic perspective on the sustainable integration of electric vehicles (EVs) and renewable resources within communities.

\bibliographystyle{elsarticle-num} 
\bibliography{RefsComBat}

\begin{thebibliography}{10}
\expandafter\ifx\csname url\endcsname\relax
  \def\url#1{\texttt{#1}}\fi
\expandafter\ifx\csname urlprefix\endcsname\relax\def\urlprefix{URL }\fi
\expandafter\ifx\csname href\endcsname\relax
  \def\href#1#2{#2} \def\path#1{#1}\fi

\bibitem{riley2023connected}
B.~Riley, L.~V. White, S.~Quilty, T.~Longden, N.~Frank-Jupurrurla,
  S.~Morton~Nabanunga, S.~Wilson, Connected: rooftop solar, prepay and reducing
  energy insecurity in remote australia, Australian Geographer (2023) 1--22.

\bibitem{li2020review}
H.~X. Li, D.~J. Edwards, M.~R. Hosseini, G.~P. Costin, A review on renewable
  energy transition in australia: An updated depiction, Journal of cleaner
  production 242 (2020) 118475.

\bibitem{malekpour2016dynamic}
A.~R. Malekpour, A.~Pahwa, A dynamic operational scheme for residential pv
  smart inverters, IEEE Transactions on Smart Grid 8~(5) (2016) 2258--2267.

\bibitem{murray2021voltage}
W.~Murray, M.~Adonis, A.~Raji, Voltage control in future electrical
  distribution networks, Renewable and Sustainable Energy Reviews 146 (2021)
  111100.

\bibitem{gholami2022fuzzy}
K.~Gholami, S.~Karimi, A.~Rastgou, Fuzzy risk-based framework for scheduling of
  energy storage systems in photovoltaic-rich networks, Journal of Energy
  Storage 52 (2022) 104902.

\bibitem{moradi2015optimal}
M.~H. Moradi, M.~Abedini, S.~R. Tousi, S.~M. Hosseinian, Optimal siting and
  sizing of renewable energy sources and charging stations simultaneously based
  on differential evolution algorithm, International Journal of Electrical
  Power \& Energy Systems 73 (2015) 1015--1024.

\bibitem{dong2020distorted}
C.~Dong, Q.~Xiao, M.~Wang, T.~Morstyn, M.~D. McCulloch, H.~Jia, Distorted
  stability space and instability triggering mechanism of ev aggregation delays
  in the secondary frequency regulation of electrical grid-electric vehicle
  system, IEEE Transactions on Smart Grid 11~(6) (2020) 5084--5098.

\bibitem{rubanenko2020analysis}
O.~Rubanenko, V.~Yanovych, Analysis of instability generation of photovoltaic
  power station, in: 2020 IEEE 7th international conference on energy smart
  systems (ESS), IEEE, 2020, pp. 128--133.

\bibitem{ali2021maximizing}
A.~Ali, K.~Mahmoud, M.~Lehtonen, Maximizing hosting capacity of uncertain
  photovoltaics by coordinated management of oltc, var sources and stochastic
  evs, International Journal of Electrical Power \& Energy Systems 127 (2021)
  106627.

\bibitem{jordehi2020energy}
A.~R. Jordehi, M.~S. Javadi, J.~P. Catal{\~a}o, Energy management in microgrids
  with battery swap stations and var compensators, Journal of Cleaner
  Production 272 (2020) 122943.

\bibitem{leveque2007investments}
F.~Lev{\^e}que, G.~Brunekreeft, Investments in generation and transmission,
  Competition and Regulation in Network Industries 2~(1) (2007) 3--8.

\bibitem{kalkbrenner2019residential}
B.~J. Kalkbrenner, Residential vs. community battery storage systems--consumer
  preferences in germany, Energy Policy 129 (2019) 1355--1363.

\bibitem{fazlhashemi2020day}
S.~S. Fazlhashemi, M.~Sedighizadeh, M.~E. Khodayar, Day-ahead energy management
  and feeder reconfiguration for microgrids with cchp and energy storage
  systems, Journal of Energy Storage 29 (2020) 101301.

\bibitem{gholami2020energy}
K.~Gholami, S.~Jazebi, Energy demand and quality management of standalone
  diesel/pv/battery microgrid using reconfiguration, International Transactions
  on Electrical Energy Systems 30~(10) (2020) e12250.

\bibitem{sheidaei2021stochastic}
F.~Sheidaei, A.~Ahmarinejad, M.~Tabrizian, M.~Babaei, A stochastic
  multi-objective optimization framework for distribution feeder
  reconfiguration in the presence of renewable energy sources and energy
  storages, Journal of Energy Storage 40 (2021) 102775.

\bibitem{yan2019robust}
X.~Yan, C.~Gu, X.~Zhang, F.~Li, Robust optimization-based energy storage
  operation for system congestion management, IEEE Systems Journal 14~(2)
  (2019) 2694--2702.

\bibitem{aryanezhad2018management}
M.~Aryanezhad, Management and coordination of ltc, svr, shunt capacitor and
  energy storage with high pv penetration in power distribution system for
  voltage regulation and power loss minimization, International Journal of
  Electrical Power \& Energy Systems 100 (2018) 178--192.

\bibitem{elkazaz2021techno}
M.~Elkazaz, M.~Sumner, E.~Naghiyev, Z.~Hua, D.~W. Thomas, Techno-economic
  sizing of a community battery to provide community energy billing and
  additional ancillary services, Sustainable Energy, Grids and Networks 26
  (2021) 100439.

\bibitem{dinh2022optimal}
N.~T. Dinh, S.~A. Pourmousavi, S.~Karimi-Arpanahi, Y.~P.~S. Kumar, M.~Guo,
  D.~Abbott, J.~A. Liisberg, Optimal sizing and scheduling of community battery
  storage within a local market, in: Proceedings of the Thirteenth ACM
  International Conference on Future Energy Systems, 2022, pp. 34--46.

\bibitem{secchi2021multi}
M.~Secchi, G.~Barchi, D.~Macii, D.~Moser, D.~Petri, Multi-objective battery
  sizing optimisation for renewable energy communities with distribution-level
  constraints: A prosumer-driven perspective, Applied Energy 297 (2021) 117171.

\bibitem{alrashidi2022community}
M.~Alrashidi, Community battery storage systems planning for voltage regulation
  in low voltage distribution systems, Applied Sciences 12~(18) (2022) 9083.

\bibitem{colyvan2008probability}
M.~Colyvan, Is probability the only coherent approach to uncertainty?, Risk
  Analysis: An International Journal 28~(3) (2008) 645--652.

\bibitem{ben2006info}
Y.~Ben-Haim, Info-gap decision theory: decisions under severe uncertainty,
  Elsevier, 2006.

\bibitem{liang2011volt}
R.-H. Liang, Y.-K. Chen, Y.-T. Chen, Volt/var control in a distribution system
  by a fuzzy optimization approach, International Journal of Electrical Power
  \& Energy Systems 33~(2) (2011) 278--287.

\bibitem{gholami2022multi}
K.~Gholami, S.~Karimi, A.~Anvari-Moghaddam, Multi-objective stochastic planning
  of electric vehicle charging stations in unbalanced distribution networks
  supported by smart photovoltaic inverters, Sustainable cities and society 84
  (2022) 104029.

\bibitem{gholami2023risk}
K.~Gholami, A.~Azizivahed, A.~Arefi, L.~Li, Risk-averse volt-var management
  scheme to coordinate distributed energy resources with demand response
  program, International Journal of Electrical Power \& Energy Systems 146
  (2023) 108761.

\bibitem{gholami2022risk}
K.~Gholami, A.~Azizivahed, A.~Arefi, Risk-oriented energy management strategy
  for electric vehicle fleets in hybrid ac-dc microgrids, Journal of Energy
  Storage 50 (2022) 104258.

\bibitem{tostado2023information}
M.~Tostado-V{\'e}liz, S.~A. Mansouri, A.~Rezaee-Jordehi, D.~Icaza-Alvarez,
  F.~Jurado, Information gap decision theory-based day-ahead scheduling of
  energy communities with collective hydrogen chain, International Journal of
  Hydrogen Energy 48~(20) (2023) 7154--7169.

\bibitem{correa2019comparative}
C.~A. Correa-Florez, A.~Michiorri, G.~Kariniotakis, Comparative analysis of
  adjustable robust optimization alternatives for the participation of
  aggregated residential prosumers in electricity markets, Energies 12~(6)
  (2019) 1019.

\bibitem{CSIRObat}
P.~Graham, J.~Hayward, J.~Foster, L.~Havas, Csiro’s gencost 2021–22,
  Consultation draft, CSIRO (2022).

\bibitem{thiruvady2019maximising}
D.~Thiruvady, C.~Blum, A.~T. Ernst, Maximising the net present value of project
  schedules using cmsa and parallel aco, in: Hybrid Metaheuristics: 11th
  International Workshop, HM 2019, Concepci{\'o}n, Chile, January 16--18, 2019,
  Proceedings 11, Springer, 2019, pp. 16--30.

\bibitem{Thiruvady2014ps}
D.~Thiruvady, M.~Wallace, H.~Gu, A.~Schutt, {A Lagrangian Relaxation and ACO
  Hybrid for Resource Constrained Project Scheduling with Discounted Cash
  Flows}, Journal of Heuristics 20~(6) (2014) 643--676.

\bibitem{combe2019cost}
M.~Combe, A.~Mahmoudi, M.~H. Haque, R.~Khezri, Cost-effective sizing of an ac
  mini-grid hybrid power system for a remote area in south australia, IET
  Generation, Transmission \& Distribution 13~(2) (2019) 277--287.

\bibitem{sun2020new}
V.~Sun, A.~Asanakham, T.~Deethayat, T.~Kiatsiriroat, A new method for
  evaluating nominal operating cell temperature (noct) of unglazed photovoltaic
  thermal module, Energy reports 6 (2020) 1029--1042.

\bibitem{dai2019igdt}
X.~Dai, Y.~Wang, S.~Yang, K.~Zhang, Igdt-based economic dispatch considering
  the uncertainty of wind and demand response, IET Renewable Power Generation
  13~(6) (2019) 856--866.

\bibitem{samimi2019robust}
A.~Samimi, N.~Rezaei, Robust optimal energy and reactive power management in
  smart distribution networks: An info-gap multi-objective approach,
  International Transactions on Electrical Energy Systems 29~(11) (2019)
  e12115.

\bibitem{gurobi}
{Gurobi Optimization, LLC}, \href{https://www.gurobi.com}{{Gurobi Optimizer
  Reference Manual}} (2023).
\newline\urlprefix\url{https://www.gurobi.com}

\end{thebibliography}

\appendix
\section{Extra Plots} \label{appendix}
This appendix presents data used to construct the models in this work. We have chosen to present all the data within this appendix, as including them would compromise the readability of the plots. Instead, a representative sample of the data, specifically for the first year of the planning horizon, is presented in Fig. \ref{fig:mylabel3}.

\begin{figure}[!h]
    \centering
    \includegraphics[width=\linewidth]{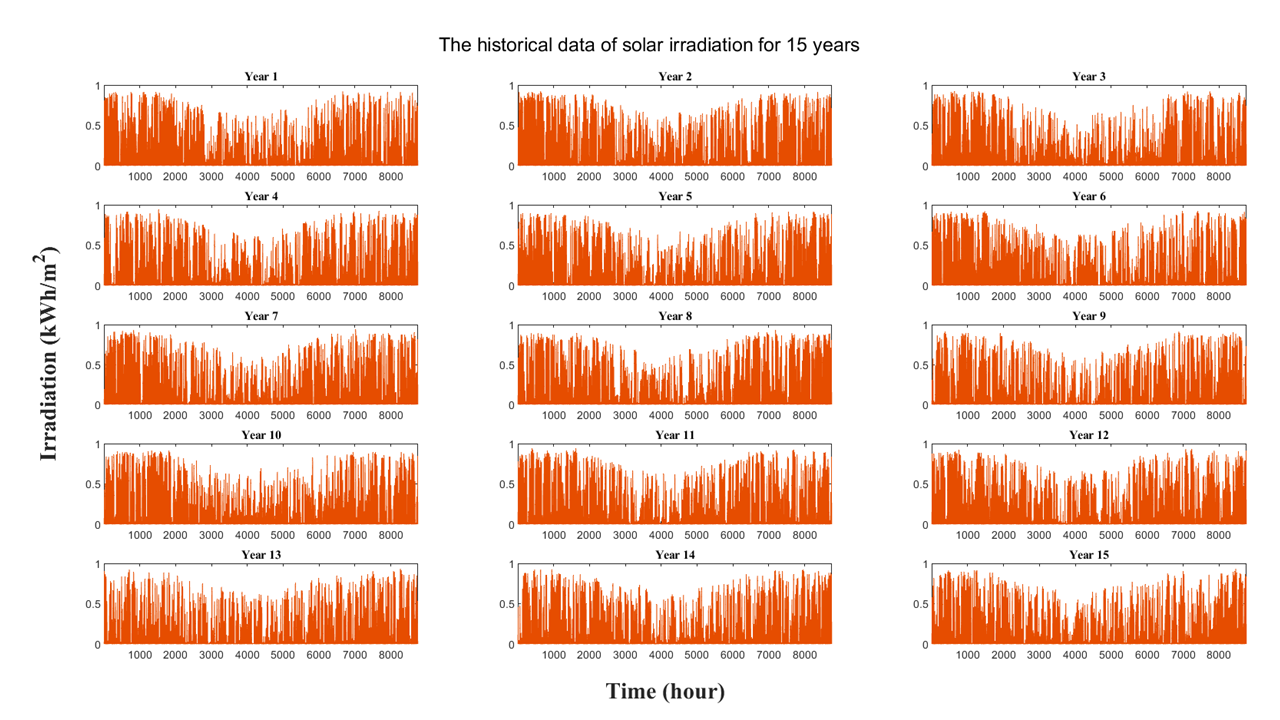}
    \caption{Historical data of solar irradiation}
    \label{fig2App}
\end{figure}

\begin{figure}[!h]
    \centering
    \includegraphics[width=\linewidth]{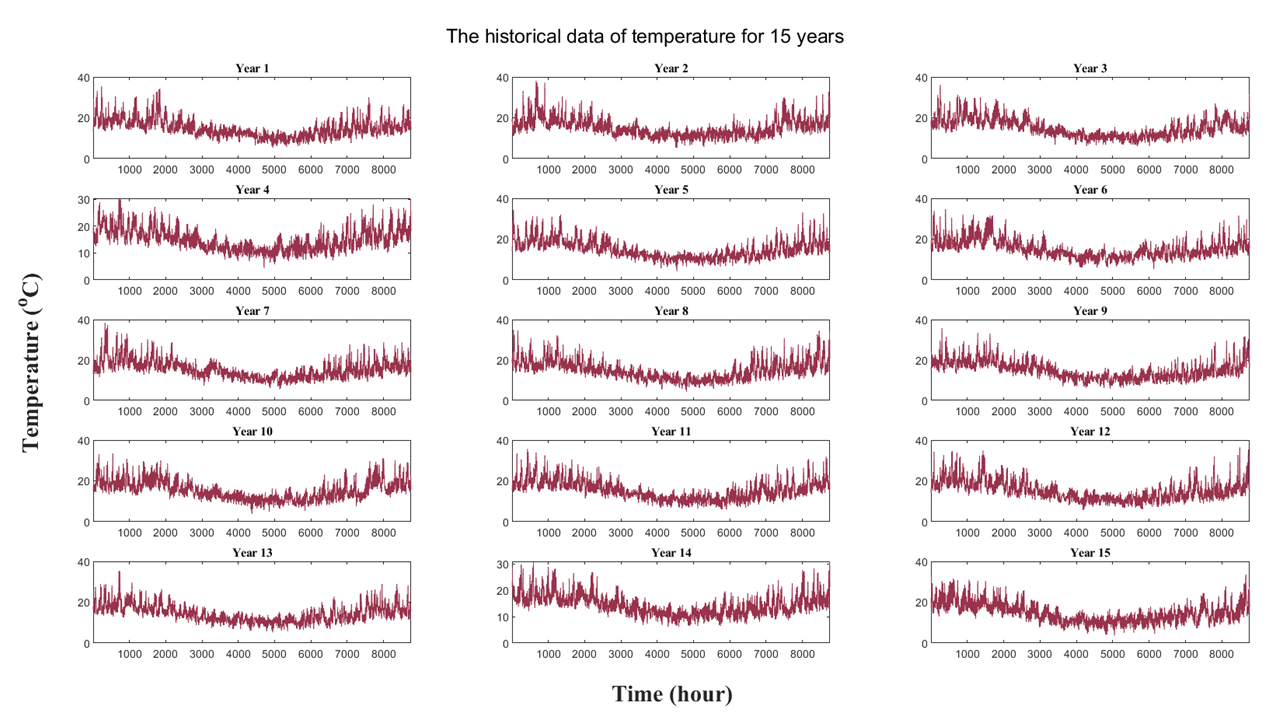}
    \caption{Historical data of temperatures}
    \label{fig3App}
\end{figure}

\begin{figure}[!h]
    \centering
    \includegraphics[width=\linewidth]{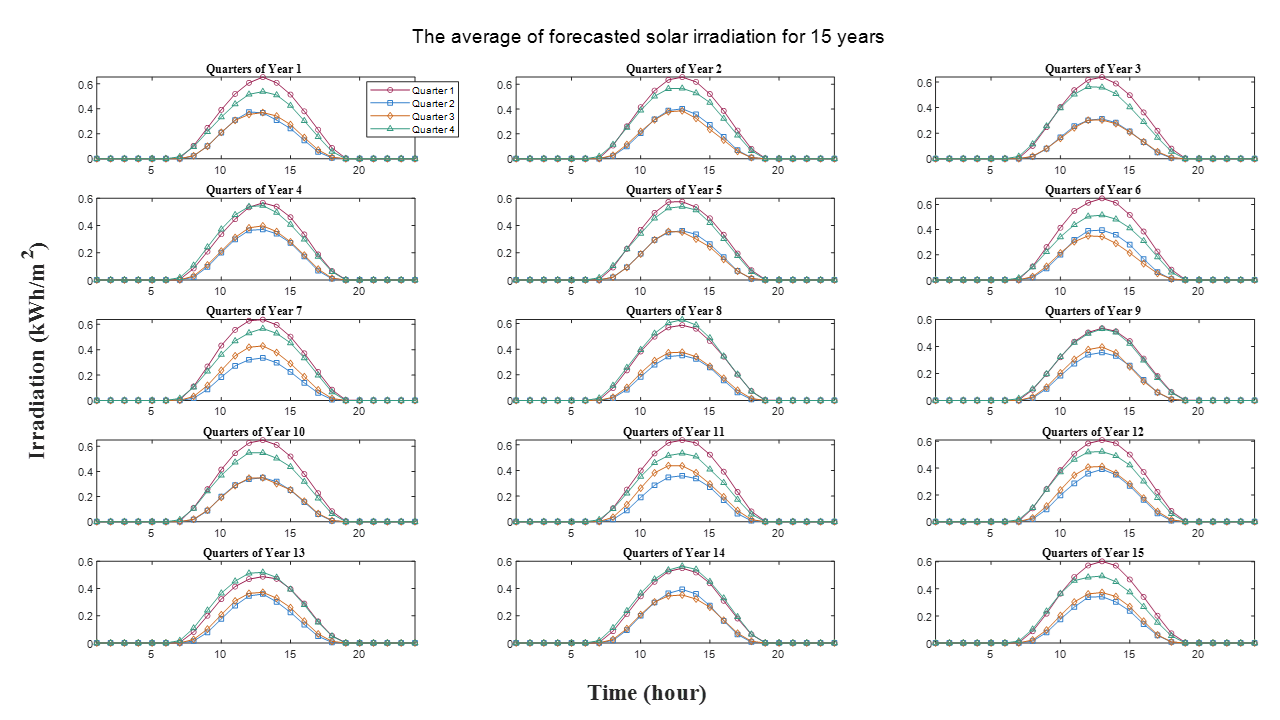}
    \caption{Quarterly Average Solar Irradiation}
    \label{fig4App}
\end{figure}

\begin{figure}[!h]
    \centering
    \includegraphics[width=\linewidth]{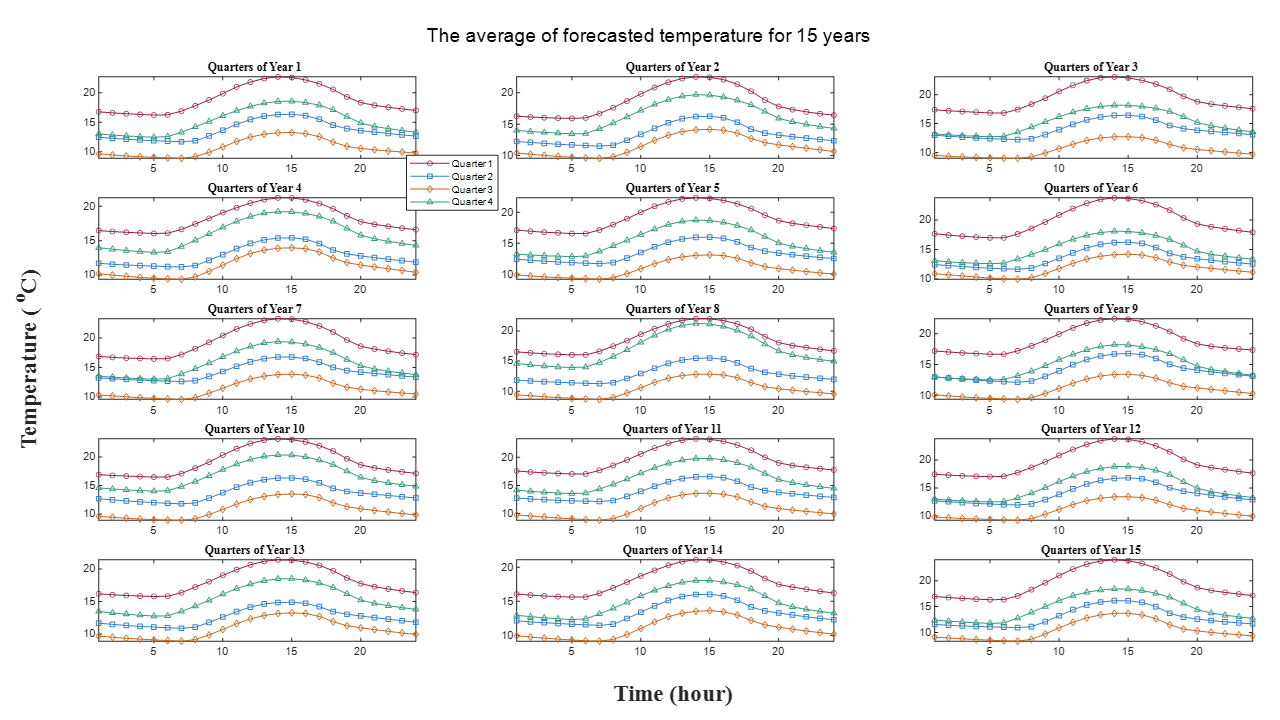}
    \caption{Quarterly average temperature}
    \label{fig5App}
\end{figure}

\begin{figure}[!h]
    \centering
    \includegraphics[width=\linewidth]{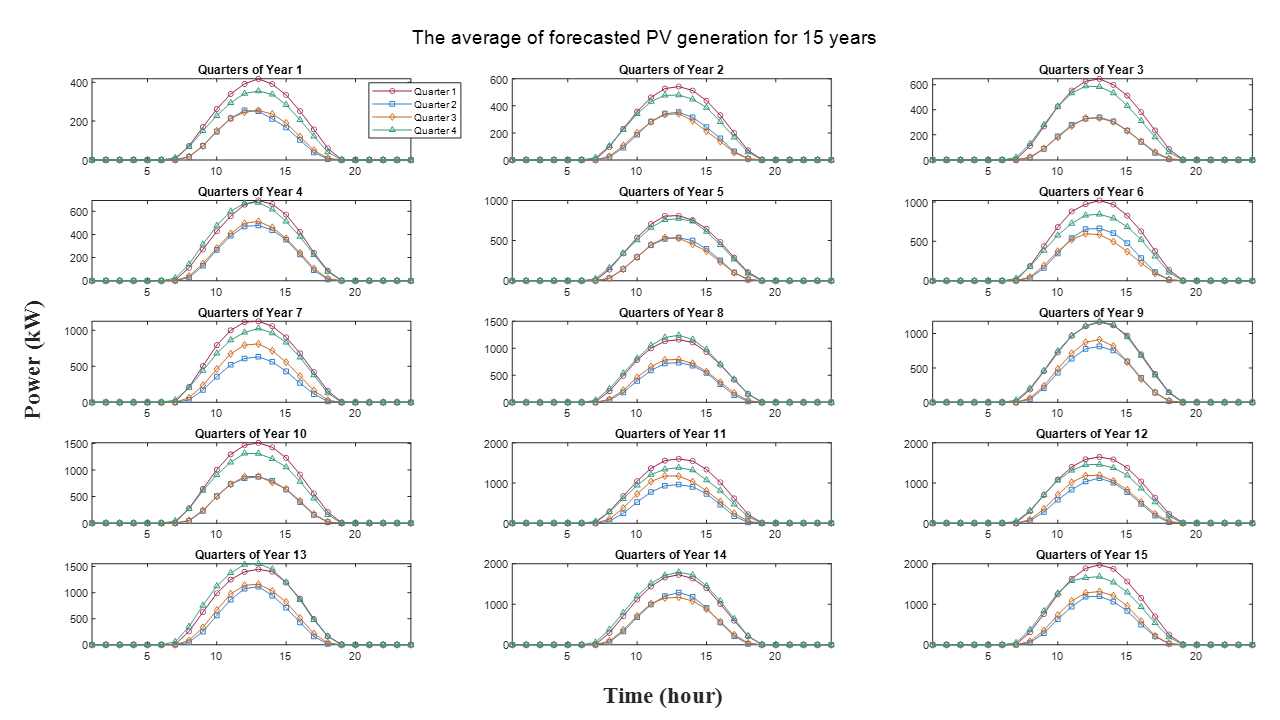}
    \caption{Quarterly PV Resource Output}
    \label{fig6App}
\end{figure}

\begin{figure}[!h]
    \centering
    \includegraphics[width=\linewidth]{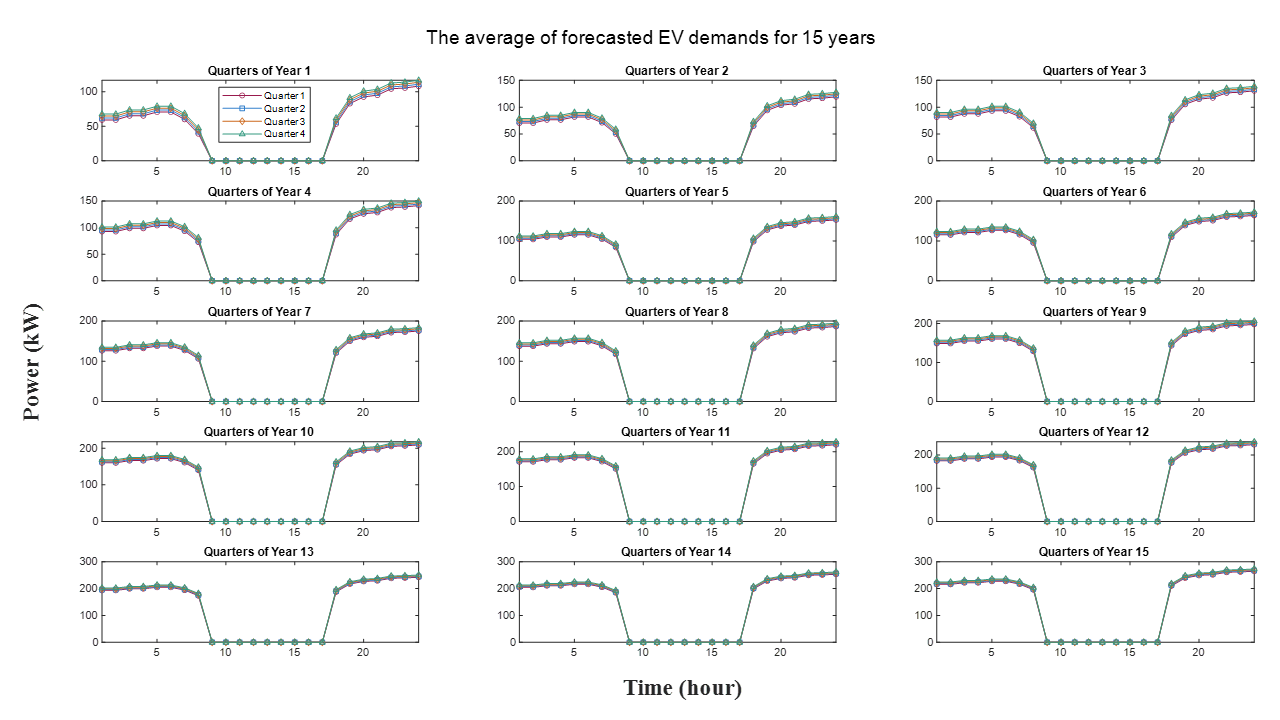}
    \caption{Quarterly Average Electric Vehicle Demands}
    \label{fig7App}
\end{figure}

\begin{figure}[!h]
    \centering
    \includegraphics[width=\linewidth]{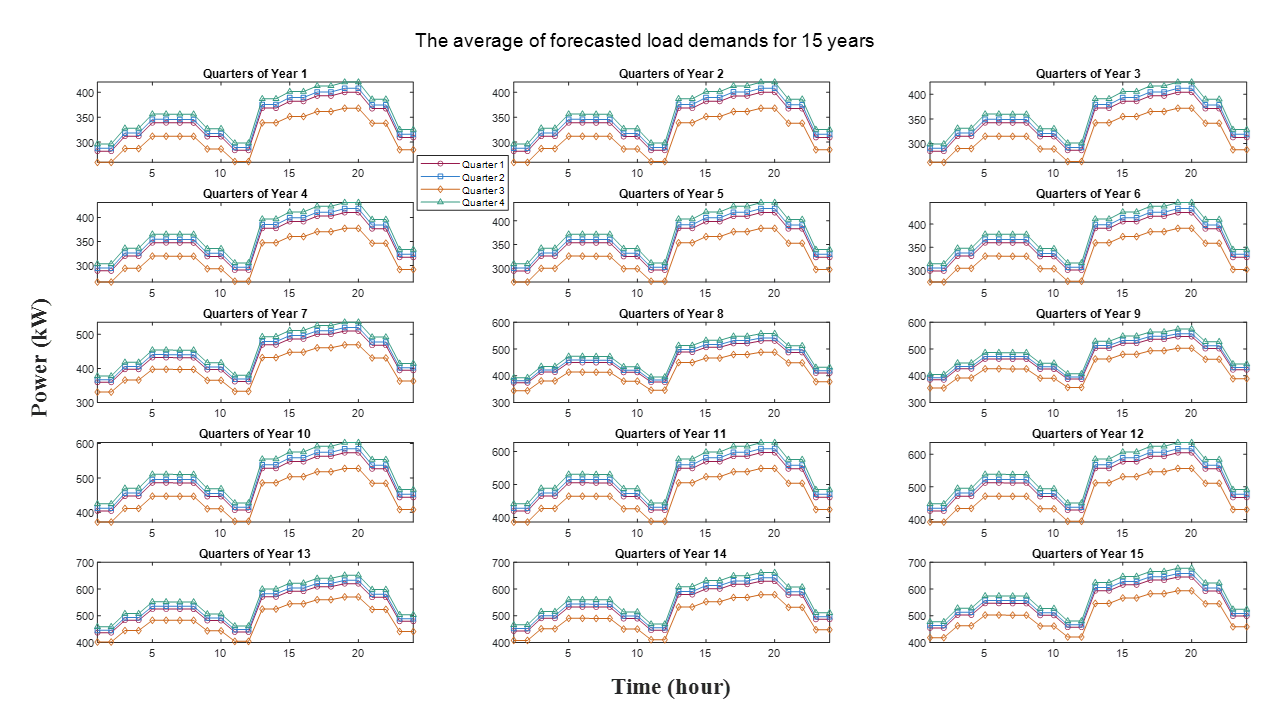}
    \caption{Quarterly average residential load demands}
    \label{fig8App}
\end{figure}

\begin{figure}[!h]
    \centering
    \includegraphics[width=\linewidth]{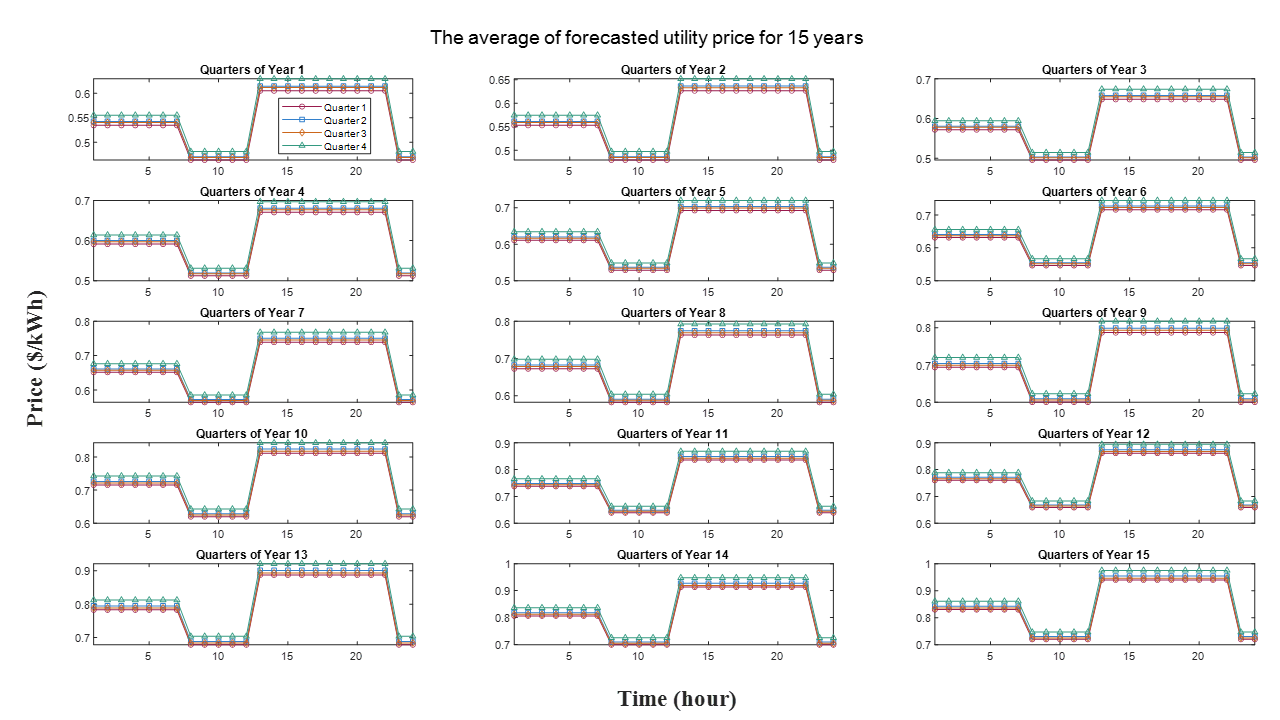}
    \caption{Quarterly Utility Prices}
    \label{fig9App}
\end{figure}

\end{document}